%% file: PtolemyDehnInvAPoly.tex
\documentclass[12pt]{amsart}
\pdfoutput=1

\usepackage{amssymb,latexsym}
\usepackage{amsmath,amsfonts,amsthm}
\usepackage{stmaryrd}
\usepackage{enumerate}
\usepackage{enumitem}
\usepackage[T1]{fontenc}
\usepackage[latin1]{inputenc}
\usepackage[english]{babel}

\usepackage{graphicx}
\usepackage{texdraw}

\usepackage[bookmarks=true,%
    colorlinks=true,%
    linkcolor=blue,%
    citecolor=blue,%
    filecolor=blue,%
    menucolor=blue,%
    urlcolor=blue,%
    breaklinks=true]{hyperref}

\usepackage{caption}

\usepackage[all]{xy}
\CompileMatrices
\def\cxymatrix#1{\xy*[c]\xybox{\xymatrix#1}\endxy}

\theoremstyle{plain}
\newtheorem{theorem}{Theorem}[section]
\theoremstyle{definition}
\newtheorem{proposition}[theorem]{Proposition}
\newtheorem{lemma}[theorem]{Lemma}
\newtheorem{definition}[theorem]{Definition}

\newtheorem{remark}[theorem]{Remark}
\newtheorem{corollary}[theorem]{Corollary}

\newtheorem{example}[theorem]{Example}

\numberwithin{equation}{section}

\newcommand{\Z}{\mathbb Z}
\newcommand{\Q}{\mathbb Q}
\newcommand{\R}{\mathbb R} 
\newcommand{\C}{\mathbb C}

\newcommand{\T}{\mathcal T}
\newcommand{\E}{\mathcal E}
\newcommand{\B}{\mathcal B}
\newcommand{\J}{\mathcal J}

\renewcommand{\P}{\mathcal P}

\DeclareMathOperator{\PSL}{PSL}
\DeclareMathOperator{\PGL}{PGL}
\DeclareMathOperator{\SL}{SL}

\DeclareMathOperator{\diag}{diag}

\DeclareMathOperator{\Conj}{Conj}

\DeclareMathOperator{\red}{red}
\DeclareMathOperator{\Stab}{Stab}
\DeclareMathOperator{\Hom}{Hom}

\usepackage{scalerel}

\usepackage{ifpdf}
\ifpdf
\fi

\usepackage[left=2.8cm,right=2.8cm,bottom=3.1cm,top=3cm]{geometry}

\title[Ptolemy coordinates, Dehn invariant, and the $A$-polynomial]{
Ptolemy coordinates, Dehn invariant and the $A$-polynomial}
\author{Christian K. Zickert}
\address{University of Maryland \\
         Department of Mathematics \\
         College Park, MD 20742-4015, USA \newline
         {\tt \url{http://www2.math.umd.edu/~zickert}}}
\email{zickert@math.umd.edu}
\thanks{
The author was supported in part by the National Science Foundation. \\
\newline
1991 {\em Mathematics Classification.} Primary 57N10, 57M27, 57M50. Secondary 13P10.
}

\date{}
\begin{document}

\begin{abstract}
We define Ptolemy coordinates for representations that are not necessarily boundary-unipotent. This gives rise to a new algorithm for computing the $\SL(2,\C)$ $A$-polynomial, and more generally the $\SL(n,\C)$ $A$-varieties. We also give a formula for the Dehn invariant of an $\SL(n,\C)$-representation. 
\end{abstract}
\maketitle

\section{Introduction}
The Ptolemy variety of a triangulated compact 3-manifold $M$ gives coordinates for boundary-unipotent representations of $\pi_1(M)$ in $\SL(n,\C)$. In this paper we define an \emph{enhanced Ptolemy variety}, which gives coordinates for representations that are not necessarily boundary-unipotent. To do this, we assume that all boundary components of $M$ are tori (this is not strictly needed; c.f.~Remarks~\ref{rm:GeneralBoundaryIntro} and \ref{rm:GeneralBoundary}), and that we have fixed a meridian and a longitude for each boundary component. The enhanced Ptolemy variety is defined via the usual Ptolemy relations, but involve additional variables $m_{i,j}$ and $l_{i,j}$, where $j=1,\dots,n-1$, and $i$ is an integer indexing the boundary components. As we shall see, the additional variables are exactly the eigenvalues of the images of the meridians and longitudes under the corresponding representation. In particular, when $n=2$, and $M$ has a single boundary torus, there are two additional variables $m$ and $l$, and by eliminating the remaining variables (e.g.~using Magma~\cite{Magma}), we obtain (a divisor of) the $\SL(2,\C)$ $A$-polynomial. In a similar fashion  one can compute (some components of) the $\SL(n,\C)$ ``$A$-varieties''. 


\subsection{The Ptolemy variety}
The Ptolemy variety $P_n(\T)$ of a $3$-manifold $M$ with a triangulation $\T$ was defined in \cite{GaroufalidisThurstonZickert} inspired by coordinates on higher Teichm\"{u}ller spaces due to Fock and Goncharov~\cite{FockGoncharov}. It is given by an ideal $I_n$ (the ideal generated by the Ptolemy relations) of a polynomial ring $Q_n$ over $\Q$, and is defined as the subset of the vanishing set of $I_n$ consisting of points where all coordinates are non-zero.

The Ptolemy variety parametrizes \emph{generically decorated} boundary-unipotent $\SL(n,\C)$-representations, where a \emph{decoration} of a boundary-unipotent representation $\rho\colon\pi_1(M)\to\SL(n,\C)$ is a $\rho$-equivariant map 
\begin{equation}\label{eq:DecorationDef}
\widehat{\widetilde M}^{(0)}\to\SL(n,\C)/N,
\end{equation}
 i.e.~an equivariant assignment of $N$-cosets (affine flags) to the vertices of each simplex of $\T$ ($N$ is the group of unipotent upper triangular matrices). A decoration is \emph{generic} if for each simplex, the four cosets assigned to the vertices are in general position. We have explicit maps
\begin{equation}\label{eq:OneToOne}
P_n(\T)\overset{1:1}{\longleftrightarrow}\left\{\txt{Generically decorated,\\boundary-unipotent\\$\pi_1(M)\to\SL(n,\C)$}\right\}\big/\Conj\to\left\{\txt{Boundary-unipotent\\$\pi_1(M)\to\SL(n,\C)$}\right\}\big/\Conj.
\end{equation}

The Ptolemy coordinates are such that the Cheeger-Chern-Simons invariant (complex volume) and Bloch invariant of a representation can be readily obtained from its Ptolemy coordinates. Explicit computations of Ptolemy varieties for the SnapPy census manifolds can be found at~\cite{Unhyperbolic}.

\begin{remark} Whether or not a representation has a generic decoration depends on the triangulation $\T$. Hence, the right map in~\eqref{eq:OneToOne} is generally not surjective. It is surjective if $\T$ is sufficiently fine (see \cite{GaroufalidisThurstonZickert}), but subdividing a triangulation creates non-ideal vertices, which increases the dimension of the preimage. If $M$ is hyperbolic and $\T$ has only essential edges, the geometric representation (and its image under the canonical irreducible representation $\SL(2,\C)\to\SL(n,\C)$) is always detected (when $n$ is even this requires obstruction classes~\cite{GaroufalidisThurstonZickert}).
\end{remark}

\subsection{The A-polynomial and its higher rank analogues}
Let $M$ be a compact $3$-manifold with boundary a torus, and let $\mu$ and $\lambda$ be a meridian and longitude.
The $A$-polynomial~\cite{APolyPaper} is (roughly) the set of eigenvalues of the images $\mu$ and $\lambda$ under the set of $\SL(2,\C)$-representations of $\pi_1(M)$. For a precise definition, consider the diagram
\begin{equation}
\cxymatrix{{X_{\SL(2,\C)}(M)\ar[d]^r&&&\\X_{\SL(2,\C)}(\partial M)&\Delta\ar[l]_-{t_\Delta}\ar[r]^-{p}_-{\cong}&\C^*\times\C^*,}}
\end{equation}
where $X_{\SL(2,\C)}$ is the character variety, $\Delta$ the set of diagonal representations of $\pi_1(\partial M)$, $r$ the restriction map, $t_\Delta$ the map sending a representation to its character, and $p$ the map taking $\rho$ to $(m,l)$ if $\rho(\mu)=\left(\begin{smallmatrix}m&\\&m^{-1}\end{smallmatrix}\right)$ and $\rho(\lambda)=\left(\begin{smallmatrix}l&\\&l^{-1}\end{smallmatrix}\right)$. 
\begin{definition}\label{def:APoly}
The \emph{$A$-polynomial} of $(M,\mu,\lambda)$ is the defining polynomial of the Zariski closure in $\C^*\times \C^*\cong \Delta$ of the union of preimages under $t_\Delta$ of the components $Y$ in $X_{\SL(2,\C)}(\partial M)$ such that $Y$ is $1$-dimensional and of the form $\overline{r(X)}$ for some component $X$ of $X_{\SL(2,\C)}(M)$.
\end{definition}
One can also define a higher rank analogue using the diagram
\begin{equation}
\cxymatrix{{X_{\SL(n,\C)}(M)\ar[d]^r&&&\\X_{\SL(n,\C)}(\partial M)&\Delta\ar[l]_-{t_\Delta}\ar[r]^-{p}_-{\cong}&(\C^*)^{n-1}\times(\C^*)^{n-1}.}}
\end{equation}
Here $p$ is the map taking $\rho$ to $(m_1,\dots, m_{n-1},l_1,\dots,l_{n-1})$ if $\rho(\mu)=\diag(m_1,\dots,m_n)$ and $\rho(\lambda)=\diag(l_1,\dots,l_n)$. 
Since every character of $\partial M$ is the character of a diagonal representation, $t_\Delta$ is surjective, and is generically $n!:1$.
\begin{definition}
The \emph{$A$-variety} of $(M,\mu,\lambda)$ is the Zariski closure in $(\C^*)^{n-1}\times(\C^*)^{n-1}$ of the union of preimages of $t_\Delta$ of the components $Y$ in $X_{\SL(n,\C)}(\partial M)$ such that $Y$ is $(n-1)$-dimensional and of the form $\overline{r(X)}$ for some component $X$ of $X_{\SL(2,\C)}(M)$.
\end{definition}
\begin{remark}
If $M$ has multiple torus boundary components, $c$ say, one can similarly define $A$-varieties by picking meridians and longitudes for each boundary component. The dimension of the $A$-variety is then $c(n-1)$.
\end{remark}

\subsection{The Dehn invariant}
The Dehn invariant was defined by Dehn~\cite{Dehn} in his solution to Hilbert's 3rd problem. In modern language (see e.g.~\cite{NeumannFest} for a historical overview) it can be viewed as a map
\begin{equation}
\nu\colon\mathcal P(\C)\to \wedge_\Z^2(\C^*),\qquad z\mapsto z\wedge(1-z),
\end{equation}
where $\mathcal P(\C)$ is the so-called \emph{pre-Bloch group}, defined as the free abelian group on $\C\setminus\{0,1\}$ subject to the \emph{five term relation}
\begin{equation}
[x]-[y]+\left[\frac{y}{x}\right]-\left[\frac{1-x^{-1}}{1-y^{-1}}\right]+\left[\frac{1-x}{1-y}\right]=0\in\mathcal P(\C).
\end{equation}
The kernel of the Dehn invariant is called the \emph{Bloch group} $\B(\C)$.
\subsection{Statement of results}
Let $M$ be an oriented, compact $3$-manifold with boundary a union of $c\geq 1$ tori together with a choice of meridian $\mu_i$ and longitude $\lambda_i$ of each boundary component. Let $\T$ be a topological ideal triangulation of $M$.
Recall that the Ptolemy variety is defined via an ideal $I_n$ generated by Ptolemy relations in a polynomial ring $Q_n$.
The enhanced Ptolemy variety $\E P_n(\T)$, defined in Section~\ref{sec:ExtendedPtolemy} below, is given by an ideal $\E I_n$ in a polynomial ring $\E Q_n$, where
\begin{equation}\label{eq:EQnMandL}
\E Q_n = \Q[\mathcal C,\mathcal M,\mathcal L],\qquad \mathcal C=\{c_{t,k}\}\qquad  \mathcal M=\{m_{i,j}\},\qquad \mathcal L=\{l_{i,j}\}.
\end{equation}
Here, $c_{t,k}$, $m_{i,j}$ and $l_{i,j}$ are free variables, i.e.~generators of $\E Q_n$. The index $i$ runs through the boundary components, and the index $j$ runs from $1$ to $n-1$. The variables $c_{t,k}$ are indexed like the standard Ptolemy coordinates. When convenient, we shall also regard these as regular functions on $\E P_n(\T)$. We shall occasionally need the functions $m_{i,n}=(m_{i,1}\cdots m_{i,n-1})^{-1}$ and $l_{i,n}=(l_{i,1}\cdots l_{i,n-1})^{-1}$.  
The standard Ptolemy variety $P_n(\T)$ embeds in $\E P_n(\T)$ as the set of points where all $m_{i,j}$'s and $l_{i,j}$'s are $1$. 

There is an action of $H^c$, where $H$ is the group of diagonal matrices, on $P_n(\T)$, which naturally extends to an action on $\E P_n(\T)$. The quotient is denoted by $\E P_n(\T)_{\red}$. A representation is \emph{boundary-Borel} if the image of each peripheral subgroup lies in a conjugate of the Borel subgroup $B\subset \SL(n,\C)$ of upper triangular matrices.
\begin{theorem}\label{thm:EPnAndReps}
We have a diagram
\begin{equation}\label{eq:OneToOneExt}
\xymatrix{\E P_n(\T)_{\red}\ar@{<->}[r]^-{1:1}&\left\{\txt{Generically decorated,\\boundary-Borel\\$\pi_1(M)\to\SL(n,\C)$}\right\}\ar[r]&\left\{\txt{Boundary-Borel\\$\pi_1(M)\to\SL(n,\C)$}\right\}}
\end{equation}
of explicit maps, which commute with~\eqref{eq:OneToOne} under the natural maps. 
\end{theorem}
 
For a manifold $N$ let $R_B(N)$ denote the set of representations of $\pi_1(N)$ in $B$ up to conjugation (in $B$). A decorated representation has a well defined peripheral holonomy (see Section~\ref{sec:PeripheralHolonomy}), and we thus have restriction maps 
\begin{equation}
r_i\colon \E P_n(\T)\to R_B(\partial_i M),
\end{equation}
which are given explicitly in terms of the Ptolemy coordinates.
Consider the maps
\begin{equation}
p_{\mathcal M,i}\colon R_B(\partial_i M)\to (\C^*)^{n-1},\qquad p_{\mathcal L,i}\colon R_B(\partial_i M)\to (\C^*)^{n-1},
\end{equation}
defined by taking $\rho\in R_B(\partial_i M)$ to $(m_{i,1},\dots,m_{i,n-1})$, respectively $(l_{i,1},\dots,l_{i,n-1})$ if
\begin{equation}
\rho_i(\mu_i)=\begin{pmatrix}m_{i,1}&*&\dots &*\\&m_{i,2}&*&\vdots\\&&\ddots&*\\&&&m_{i,n}\end{pmatrix},\qquad \rho_i(\lambda_i)=\begin{pmatrix}l_{i,1}&*&\dots &*\\&l_{i,2}&*&\vdots\\&&\ddots&*\\&&&l_{i,n}\end{pmatrix}.
\end{equation}
\begin{theorem}\label{thm:Projection}
We have a commutative diagram
\begin{equation}
\cxymatrix{{\E P_n(\T)\ar[rd]^{\Pi_{\mathcal M,\mathcal L}}\ar[d]_-r&\\R_B(\partial_1 M)\times\dots\times R_B(\partial_h M)\ar[r]^-p&(\C^*)^{2h(n-1)}}}
\end{equation}
where $r=(r_1,\dots, r_c)$, $p=(p_{1,\mathcal M},\dots,p_{c,\mathcal M},p_{1,\mathcal L},\dots,p_{c,\mathcal L})$,
and $\Pi_{\mathcal M,\mathcal L}$ is projection onto the $\mathcal M$ and $\mathcal L$ coordinates.\qed
\end{theorem}
\begin{corollary}
If $M$ has a single torus boundary, each $1$-dimensional component of the Zariski closure of the image of $\E P_2(\T)$ in $\C^*\times\C^*$ is a component of the $\SL(2,\C)$ $A$-polynomial.\qed
\end{corollary}

\begin{remark} It follows from~\cite[Theorem~1.8]{GaroufalidisThurstonZickert} that if $M$ is hyperbolic, the geometric component of the $A$-polynomial is always detected.
\end{remark}

\begin{remark}
The components of the image of $\E P_n(\T)$ in $(\C^*)^{2c(n-1)}$ contain slightly finer information than the corresponding component of the $A$-variety. The $A$-variety is by definition symmetric under the Weil group action, but this need not be the case here. For example,
there may exist a decorated representation $\rho$ whose image $\rho_i$ in $R_B(\partial_i M)$ satisfies 
\begin{equation}
\rho_i(\mu_i)\begin{pmatrix}m&1&\\&m&1\\&&m\end{pmatrix},\qquad \rho_i(\lambda_i)=\begin{pmatrix}l_1&&\\&l_2&\\&&l_3\end{pmatrix},
\end{equation}
but there may not exist a representation $\tau$ with $\tau_i(\lambda_i)=\diag(l_3,l_2,l_1)$. 
\end{remark}

\subsubsection{Formula for the Dehn invariant}
In~\cite{GaroufalidisThurstonZickert} we defined a map
\begin{equation}\label{eq:PtolemyToBloch}
P_n(\T)\to \B(\C).
\end{equation}
assigning a Bloch group element to a decorated representation. 

\begin{theorem}\label{thm:DehnInvariantFormula}
There is a map 
\begin{equation}
\lambda\colon\E P_n(\T)\to\mathcal P(\C)
\end{equation}
extending the map~\eqref{eq:PtolemyToBloch}. For each $x\in\E P_n(\T)$, we have 
\begin{equation}\label{eq:DehnInvariantFormula}
\nu\circ\lambda(x)=-\sum_{i=1}^c\sum_{j,k=1}^{n-1}A^{-1}_{jk}\overline m_{i,j}(x)\wedge \overline l_{i,k}(x)\in\wedge^2(\C^*),
\end{equation}
where $A$ is the Cartan matrix of $\SL(n,\C)$, $\overline m_{i,j}=m_{i,j}/m_{i,j+1}$, and $\overline l_{i,j}=l_{i,j}/l_{i,j+1}$. In particular, the Dehn invariant of a decorated representation depends only on the restriction to the boundary. 
\end{theorem}
\begin{remark}
The formula is well defined since $\wedge^2_\Z(\C^*)$ is a $\Q$-vector space. The $jk$'th entry of $A^{-1}$ is $j(n-k)/n$ when $j\leq k$.
\end{remark}

\begin{remark}
An analogue of~\eqref{eq:DehnInvariantFormula} was proved for $n=3$ in~\cite{BergeronFalbelGuilloux}, and the general formula was stated without proof in~\cite[(3.50)]{DimofteGabellaGoncharov}. The proof is an elementary consequence of the symplectic properties of the gluing equations proved in~\cite{GaroufalidisZickert}, and also~\cite{Guilloux}. 
\end{remark}

\begin{remark}
Much of this theory could also be stated for $\PGL(n,\C)$-representations using the gluing equation varieties $V_n(\T)$ defined in~\cite{GaroufalidisGoernerZickert}. This would generalize the relationship between Thurston's gluing equations and the $\PSL(2,\C)$ $A$-polynomial (see e.g.~\cite{Champanerkar}). The enhanced Ptolemy variety further illustrates the duality between shapes and Ptolemy coordinates discussed in~\cite{GaroufalidisGoernerZickert}; boundary-unipotent $\PGL(n,\C)$-representations are parametrized by introducing more \emph{equations} (cusp equations), whereas non-boundary-unipotent $\SL(n,\C)$-representations are parametrized by introducing more \emph{variables}.
\end{remark}

\begin{remark}\label{rm:GeneralBoundaryIntro} The condition that all boundary components be tori is not strictly needed. We have restricted to this case for ease of exposition, and because this is the case of most general interest. See Remark~\ref{rm:GeneralBoundary} for a description of how the theory is modified in the general case.
\end{remark}

\subsection{Acknowledgements} The author wishes to thank Stavros Garoufalidis, Matthias Goerner and Dylan Thurston for useful conversations, and Fabrice Rouillier for assistance with some of the computations.

\section{The (standard) Ptolemy variety}\label{sec:StandardPtolemy}
To set up our notation, we first recall the definition of the standard Ptolemy variety $P_n(\T)$ \cite{GaroufalidisThurstonZickert,GaroufalidisGoernerZickert}. We refer to~\cite{GaroufalidisGoernerZickert} for conventions regarding triangulations.

\subsection{Preliminaries}
Let $M$ be a compact, oriented $3$-manifold, and let $\T$ be a topological ideal triangulation, i.e.~a triangulation where the $0$-cells correspond to boundary components (if $\partial M=\emptyset$, a triangulation of $M$ may be regarded as an ideal triangulation of a manifold with boundary a union of spheres).

Identify each simplex of $\T$ with a standard simplex
\begin{equation}
\Delta^3_n=\big\{(x_0,x_1,x_2,x_3)\subset \R^4\bigm\vert 0\leq x_i\leq n,\enspace x_0+x_1+x_2+x_3=n \big\}.
\end{equation}
Let $\Delta^3_n(\Z)$ be the \emph{integral points} of $\Delta^3_n$, and let $\dot\Delta^3_n(\Z)$ be the integral points with the $4$ vertex points removed.
An \emph{integral point} of $M$ is an equivalence class under face pairings of the (non-vertex) integral points of the simplices of $M$.

\subsection{The Ptolemy relations}
Assign to each integral point $t\in\dot\Delta^3_n(\Z)$ of each simplex $\Delta_k$ of $\T$ a \emph{Ptolemy coordinate} $c_{t,k}$. We shall regard the $c_{t,k}$'s as generators of a polynomial ring $Q_n$ over $\Q$. 
\begin{definition}
The \emph{Ptolemy relations} are the relations
\begin{equation}\label{eq:PtolemyRelations}
c_{s+1001,k}c_{s+0110,k}-c_{s+1010,k}c_{s+0101,k}+c_{s+1100,k}c_{s+0011,k}=0,\quad s\in \Delta_{n-2}^3, \enspace\Delta_k\in\T.
\end{equation}
\end{definition}
\begin{remark} When $n=2$, we denote the Ptolemy coordinates by $c_{ij,k}$ instead of $c_{t,k}$ (e.g.~$c_{02,k}$ instead of $c_{1010,k}$) to simplify the notation. 
\end{remark}
\subsection{The identification relations}
Recall that $M$ is obtained by pairing faces of the simplices of $\T$. Each face pairing may be encoded by a permutation $\sigma\in S_4$. The natural identification of $S_4$ with the symmetry group of an ordered simplex induces an action of $S_4$ on $\Delta^3_n(\Z)$.
 
For $t=(t_0,\dots,t_3)\in\Delta^3_n(\Z)$ we can write the $n\times n$ identity matrix $I$ as a concatenation $I=[I^t_0|I^t_1|I^t_2|I^t_3]$ of four (possibly empty) matrices $I^t_i$ of respective size $n\times t_i$. For $\sigma\in S_4$ define 
\begin{equation}
I_{\sigma,t}=[I^t_{\sigma(0)}|I^t_{\sigma(1)}|I^t_{\sigma(2)}|I^t_{\sigma(3)}].
\end{equation} 
Clearly, the determinant of $I_{\sigma,t}$ is either $1$ or $-1$.
\begin{definition} 
An \emph{identification relation} is a relation of the form
\begin{equation}\label{eq:IdentificationRelations}
c_{t_j,j}=\det(I_{\sigma,t_k})c_{t_k,k}
\end{equation}
defined whenever two integral points $t_j$ and $t_k$ of $\Delta_j$ and $\Delta_k$, respectively, are identified in $M$ via a face pairing with permutation $\sigma\in S_4$.
The determinant $\det(I_{\sigma,t_k})$ is called a \emph{sign multiplier}.
\end{definition}

\begin{figure}[htb]
\begin{minipage}[b]{0.48\textwidth}
\input{figures_gen/IdentificationRelations_n=3.tex}
\end{minipage}
\begin{minipage}[b]{0.48\textwidth}
\input{figures_gen/IdentificationRelations_n=4.tex}
\end{minipage}
\\
\begin{minipage}[t]{0.5\textwidth}
\captionsetup{width=0.8\textwidth}
\caption{Identification relations for $n=3$, e.g.~$c_{1011,0}=-c_{0111,1}$.}\label{fig:IdRel3}
\end{minipage}
\begin{minipage}[t]{0.48\textwidth}
\captionsetup{width=0.9\textwidth}
\caption{Identification relations for $n=4$, e.g.~$c_{0031,0}=-c_{0103,1}$.}\label{fig:IdRel4}
\end{minipage}
\end{figure}

The sign multipliers are illustrated in Figures~\ref{fig:IdRel3} and~\ref{fig:IdRel4}. For $n=2$, the sign multiplier is positive and the relation is $c_{ij,k}=c_{i^\prime j^\prime,k^\prime}$, if and only if $i-j$ and $i^\prime-j^\prime$ have the same sign. The sign multipliers are always positive (for any $n$) for an \emph{ordered} triangulation.

\begin{definition}
The Ptolemy relations~\eqref{eq:PtolemyRelations} and the identification relations~\eqref{eq:IdentificationRelations} define an ideal $I_n$ in a polynomial ring $Q_n$. The \emph{Ptolemy variety} is the Zariski open subset of the vanishing set of $I_n$ consisting of the points where all Ptolemy coordinates are non-zero.
\end{definition}

\begin{remark} By selecting once and for all a representative for each integral point of $M$, one can eliminate the identification relations. The Ptolemy variety is thus given by an ideal generated by Ptolemy relations in a polynomial ring with a variable for each integral point of $M$.
\end{remark}

\subsection{The diagonal action}\label{sec:DiagonalAction}
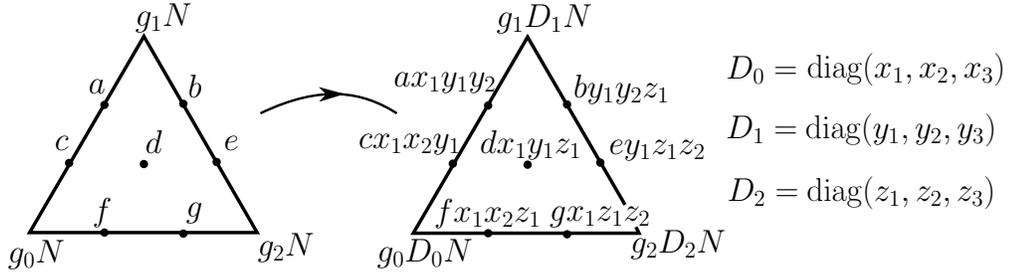
\begin{figure}[htb]
\hspace{-1cm}\scalebox{0.6}{\input{figures_gen/DiagonalAction.tex}}
\caption{The diagonal action on decorations and Ptolemy coordinates.}\label{fig:DiagonalAction}
\end{figure}

Let $H$ denote the set of diagonal matrices in $\SL(n,\C)$ and let $c$ denote the number of boundary components of $M$. It follows from~\eqref{eq:DecorationDef} that $H^c$ acts on the set of decorations of a boundary-unipotent representation $\rho$ by right multiplication. Hence, by~\eqref{eq:OneToOne}, $H^c$ also acts on $P_n(\T)$. The action is given by
\begin{equation}\label{eq:DiagonalAction}
c_{t,k}\mapsto\det\big(\{D_{k,0}\}_{t_0}\cup\{D_{k,1}\}_{t_1}\cup\{D_{k,2}\}_{t_2}\cup\{D_{k,3}\}_{t_3}\big)c_{t,k},
\end{equation}
where $D_{k,i}$ denotes the diagonal matrix of the $i$th vertex of $\Delta_k$, and $\{A\}_k$ is the ordered set consisting of the first $k$ column vectors of $A$. The action is illustrated in Figure~\ref{fig:DiagonalAction}.

The quotient $P_n(\T)_{\red}$ is called the \emph{reduced Ptolemy variety}. It is shown in~\cite{PtolemyField} that $P_n(\T)_{\red}$ is given by an ideal obtained from $I_n$ by adding $c(n-1)$ relations of the form $c_{t,k}=1$ for suitably chosen Ptolemy coordinates.
Since we usually only care about the representation and not the decoration, this simplifies the computations.

\subsection{Example, $m004$}
For the triangulation of $m004$ (the figure $8$ knot) given in Figure~\ref{fig:m004Triangulation} and $n=3$, there are $8$ integral points ($4$ edge points and $4$ face points). The identification relations are
\begin{equation}
\begin{aligned}
&c_{2100,1}=c_{1002,0}=c_{2010,1}=c_{1020,0}=c_{0102,1}=c_{0120,0},\qquad&c_{1101,1}&=c_{1011,0}\\
&c_{1200,1}=c_{2001,0}=c_{1020,1}=c_{2010,0}=c_{0201,1}=c_{0210,0},\qquad&c_{1011,1}&=c_{1101,0}\\
&c_{0120,1}=c_{2100,0}=c_{0021,1}=c_{0012,0}=c_{2001,1}=c_{0102,0},\qquad&c_{1110,1}&=-c_{1110,0}\\
&c_{0210,1}=c_{1200,0}=c_{0012,1}=c_{0021,0}=c_{1002,1}=c_{0201,0},\qquad&c_{0111,1}&=-c_{0111,0}
\end{aligned}
\end{equation}
and the Ptolemy relations are
\begin{equation}
\begin{gathered}
c_{2001,k}c_{1110,k}-c_{2010,k}c_{1101,k}+c_{2100,k}c_{1011,k}=0\\
c_{1101,k}c_{0210,k}-c_{1110,k}c_{0201,k}+c_{1200,k}c_{0111,k}=0\\
c_{1011,k}c_{0120,k}-c_{1020,k}c_{1011,k}+c_{1110,k}c_{0021,k}=0\\
c_{1002,k}c_{0111,k}-c_{1011,k}c_{0102,k}+c_{1101,k}c_{0012,k}=0
\end{gathered}
\end{equation}

Eliminating the identification relations, we obtain a system 
of $8$ equations in $8$ variables (one for each integral point of $M$). The reduced Ptolemy variety $P_3(\T)_{\red}$ (see Section~\ref{sec:DiagonalAction}) is obtained by adding two additional relations (e.g.~$c_{0012,0}=1$ and $c_{0111,0}=1$). It is zero dimensional and has $3$ components defined over $\Q(\sqrt{-7})$, $\Q(\sqrt{-7})$ and $\Q(\sqrt{-3})$, respectively. The components can be computed in SnapPy (with Sage of Magma installed) or retrieved from the data base in~\cite{Unhyperbolic}. We shall not display them here.

\section{The enhanced Ptolemy variety}\label{sec:ExtendedPtolemy}
We now assume that each boundary component $\partial_i M$ of $M$ is a torus, and that we have fixed once and for all meridians $\mu_i$ and longitudes $\lambda_i$ of $\partial_i M$.
Assign Ptolemy coordinates $c_{t,i}$ to the integral points of the simplices of $\T$ as in section~\ref{sec:StandardPtolemy}, and let 
\begin{equation}
\E Q_n=\Q[\mathcal C,\mathcal M,\mathcal L]=Q_n\otimes \Q[\mathcal M,\mathcal L] 
\end{equation}
with $\mathcal C$, $\mathcal M$ and $\mathcal L$ as in~\eqref{eq:EQnMandL}. The Ptolemy relations are defined as in~\eqref{eq:PtolemyRelations}, but the identification relations now involve the additional variables $m_{i,j}$ and $l_{i,j}$.

Choose a fundamental rectangle $R_i$ in $\widehat M$ of each boundary component. The triangulation of $M$ by truncated simplices induces a triangulation of each $R_i$ (see Figures~\ref{fig:FacePairings} and~\ref{fig:FundamentalRectangles}). It will later become apparent that a change of fundamental rectangles changes the variety by a canonical isomorphism.

\begin{figure}[htb]
\begin{center}
\hspace{-4mm}
\begin{minipage}[t]{0.48\textwidth}
\scalebox{0.7}{\input{figures_gen/FacePairing.tex}}
\end{minipage}
\hspace{2.5cm}
\begin{minipage}[b]{0.35\textwidth}
\scalebox{0.5}{\input{figures_gen/CuspImages.tex}}
\end{minipage}
\\
\begin{minipage}[t]{0.48\textwidth}
\caption{Face pairings.}\label{fig:FacePairings}
\end{minipage}
\begin{minipage}[t]{0.48\textwidth}
\caption[Bla]{Fundamental rectangles.}\label{fig:FundamentalRectangles}
\end{minipage}
\end{center}
\end{figure}

The sides $\mu_i^\prime$ and $\lambda_i^\prime$ of the rectangles $R_i$ map to generators of $H_1(\partial_i M)$, so there exist unique integers $a_i$, $b_i$, $c_i$ and $d_i$ so that
\begin{equation}\label{eq:MuMuPrime}
\mu_i^\prime=\mu_i^{a_i}\lambda_i^{b_i}, \qquad 
\lambda_i^\prime=\mu_i^{c_i}\lambda_i^{d_i},\qquad \det\begin{pmatrix}a_i&b_i\\c_i&d_i\end{pmatrix}=\pm 1.
\end{equation} 
Define auxilary variables $m_{i,j}^\prime$ and $l_{i,j}^\prime$ by
\begin{equation}
m_{i,j}^\prime=m_{i,j}^{a_i}l_{i,j}^{b_i},\qquad l_{i,j}^\prime=m_{i,j}^{a_i}l_{i,j}^{c_i}.
\end{equation}
They are rational functions in $m_{i,j}$ and $l_{i,j}$.
Also, let 
\begin{equation}
D_{\mu_i^\prime}=\diag(m_{i,1}^\prime,\dots,m_{i,n}^\prime),\qquad D_{\lambda_i^\prime}=\diag(l_{i,1}^\prime,\dots,l_{i,n}^\prime)
\end{equation}

Each face pairing $\alpha$ pairs a face $f_j$ of a simplex $\Delta_j$ with a face $f_i$ of a simplex $\Delta_i$ (possibly $i=j$). We shall assign to $\alpha$ a tuple $M_\alpha=(M_{\alpha0},M_{\alpha1},M_{\alpha2},M_{\alpha3})$ of $n\times n$ matrices, which will determine how the Ptolemy coordinates on $f_j$ and $f_i$ are identified. We define $M_{\alpha i}$ as follows: If $v$ is the vertex, not in $f_j$, let $M_{\alpha v}=I$. For example, for $a$ and $b$ in Figure we have $M_{a3}=M_{b2}=I$. Each vertex in $f_j$ corresponds to a triangle in one of the $R_i$'s, and the face pairing pairs this triangle with an adjacent triangle. If for a vertex $v$ in $f_j$, the face pairing is inside $R_k$, let $M_{\alpha v}=I$. Otherwise let $M_{\alpha v}$ be either $D_{\mu_k^\prime}$ or $D_{\lambda_k^\prime}$, or their respective inverses according to whether the face pairing is via $\mu_k^\prime$, $\lambda_k^\prime$, or their inverses. For example, we have
\begin{equation}
M_a=(I,D_{\mu_k^\prime},I,I),\qquad M_b=(I,I,I,D_{\lambda_l^\prime}).
\end{equation} 

\subsection{Orientation conventions}
We orient each boundary component with the counter-clockwise orientation as viewed from the ideal point, and orient the fundamental rectangles accordingly. In particular, we have $\langle \mu_i^\prime,\lambda_i^\prime\rangle=1$, where $\langle,\rangle$ is the intersection form.
For knot complements it is more common to orient the boundary with the counter-clockwise orientation as viewed from within the manifold, so in this case $\langle\mu,\lambda\rangle=-1$ for the standard meridian and longitude.
Unless otherwise specified, we assume that $\mu_i$ and $\lambda_i$ are chosen such that $\langle \mu_i,\lambda_i\rangle=1$, so that the determinants in~\eqref{eq:MuMuPrime} are $1$.

\subsection{Identification relations}\label{sec:ExtendedIdentification}

For a matrix $g$ and a non-negative integer $i$, let $\{g\}_i$ be the ordered set consisting of the first $i$ column vectors of $g$. 
\begin{definition}
The \emph{peripheral multiplier} of a face pairing $\alpha$ is the map
\begin{equation}
c(M_\alpha)\colon \Delta^3_n(\Z)\to \E Q_n,\qquad t=(t_0,t_1,t_2,t_3)\mapsto c(M_\alpha)_t,
\end{equation}
where
\begin{equation}\label{eq:PtolemyMAlpha}
c(M_\alpha)_t=\det\big(\{M_{\alpha0}\}_{t_0}\cup\{M_{\alpha1}\}_{t_1}\cup \{M_{\alpha2}\}_{t_2}\cup \{M_{\alpha3}\}_{t_3}\big)
\end{equation}
\end{definition}
The peripheral multiplier is illustrated in Figures~\ref{fig:Multiplier2} and~\ref{fig:Multiplier3}. It is simply the Ptolemy assignment (see~\eqref{eq:DecToPtolemy} c.f.~\cite{GaroufalidisThurstonZickert,GaroufalidisGoernerZickert}) associated to the tuple $M_\alpha$.

\begin{definition}
An \emph{identification relation} is a relation of the form
\begin{equation}\label{eq:ExtendedIdentifications}
c_{t_i,i}=\det(I_{\sigma,t_j})c(M_\alpha)_{t_j}c_{t_j,j}
\end{equation}
defined whenever two integral points $t_i$ and $t_j$ of $\Delta_i$ and $\Delta_j$, respectively, are identified in $M$ via a face pairing $\alpha$ with face pairing permutation $\sigma\in S_4$.
\end{definition}

\begin{figure}[htb]
\begin{minipage}[b]{0.48\textwidth}
\input{figures_gen/ExtIdenRel_n=2.tex}
\end{minipage}
\begin{minipage}[b]{0.48\textwidth}
\input{figures_gen/ExtIdenRel_n=3.tex}
\end{minipage}
\\
\begin{minipage}[t]{0.5\textwidth}
\captionsetup{width=0.8\textwidth}
\caption{Peripheral multiplier for $n=2$.}\label{fig:Multiplier2}
\end{minipage}
\begin{minipage}[t]{0.48\textwidth}
\captionsetup{width=0.9\textwidth}
\caption{Peripheral multiplier for $n=3$.}\label{fig:Multiplier3}
\end{minipage}
\end{figure}

\begin{definition}
The \emph{enhanced Ptolemy ideal} is the ideal $\E I_n$ generated by the Ptolemy relations~\eqref{eq:PtolemyRelations} and the identification relations~\eqref{eq:ExtendedIdentifications}.
The Ptolemy variety is the subset of the vanishing set of $\E I_n$ consisting of points where all coordinates are non-zero.
\end{definition}

\begin{remark}\label{rm:GeneralBoundary}
If the $i$th boundary components has genus $g_i$, one can still define the enhanced Ptolemy variety after picking generators $\mu_{i,1},\dots,\mu_{i,g_i}$ and $\lambda_{i,1},\dots,\lambda_{i,g_i}$ for $H_1(\partial_i M)$, and fundamental domains for $\widetilde{\partial M_i}$ (if $g_i$=0, no choice is necessary). The tuples $M_a$ used to define the identification relations can still be defined purely from the combinatorics of the triangulated fundamental domains. We restrict to the case of torus boundary components to avoid excessive indexing, and because this case is the most interesting.
\end{remark}

\subsection{The diagonal action}
The diagonal action~\eqref{eq:DiagonalAction} of $H^c$ on $P_n(\T)$ extends to an action on $\E P_n(\T)$ given by the same formula. A comparison of~\eqref{eq:DiagonalAction} and~\eqref{eq:PtolemyMAlpha} shows that the action respects the identification relations, and is thus well defined.

\subsection{Dehn invariant and gluing equations}\label{sec:DehnInvariant}
In~\cite{GaroufalidisGoernerZickert} we defined a variety $V_n(\T)$ generalizing Thurston's gluing equation variety, and a monomial map \begin{equation}
\mu\colon P_n(\T)\to V_n(\T),\qquad \mu(c)_s^{1100}=\frac{c_{s+1001}c_{s+0110}}{c_{s+1010}c_{s+0101}},
\end{equation}
given by assigning a shape to each subsimplex (see Section~\ref{sec:ShapeAssignments} for a review). By the same formula, we obtain a map $\mu\colon\E P_n(\T)\to V_n(\T)$. The fact that the shape assignment $\mu(c)$ satisfies the generalized gluing equations~\eqref{eq:GeneralizedGluingEquations} follows from the same cancelation argument as in~\cite{GaroufalidisGoernerZickert}. There is a map
\begin{equation}
\lambda\colon V_n(\T)\to \P(\C),\qquad z\mapsto \sum z_{s,\Delta}^{1100},
\end{equation}
taking a shape assignment to the formal sum of the shapes of the subsimplices.
A shape assignment thus has an associated element in $\P(\C)$, and by composing with $\mu$, the same holds for a Ptolemy assignment.
\begin{definition}
The \emph{Dehn invariant} of a shape assignment $z\in V_n(\T)$ is the Dehn invariant of its associated element in $\P(\C)$. The Dehn invariant of a Ptolemy assignment $c\in\E P_n(\T)$ is defined similarly.
\end{definition}

\subsection{Extended Ptolemy variety of the figure $8$ knot, $n=2$}
Consider the triangulation of $m004$ given in Figure~\ref{fig:m004Triangulation} and the fundamental rectangle given in Figure~\ref{fig:m004Rectangle}. The knot meridian $\mu$ and longitude $\lambda$ are given, respectively, by $\mu^{\prime}=\mu$, and $\lambda^{\prime}=\mu^2\lambda^{-1}$.
\begin{figure}[htb]
\hspace{-1cm}
\begin{minipage}[b]{0.49\textwidth}
\scalebox{0.55}{\input{figures_gen/m004Triangulation.tex}}
\end{minipage}
\hspace{2cm}
\begin{minipage}[t]{0.4\textwidth}
\scalebox{0.47}{\input{figures_gen/m004Cusp.tex}}
\end{minipage}
\\
\begin{minipage}[t]{0.48\textwidth}
\caption{Triangulation of $m004$.}\label{fig:m004Triangulation}
\end{minipage}
\begin{minipage}[t]{0.48\textwidth}
\caption{Fundamental rectangle.}\label{fig:m004Rectangle}
\end{minipage}
\end{figure}

We have
\begin{equation}
M_a=(D_{\mu^\prime},D_{\lambda^\prime}^{-1},I,D_{\mu^\prime}),\qquad M_b=M_d=(I,I,I,I),\qquad M_c=(I,D_{\mu^\prime},D_{\mu^\prime},I).
\end{equation}

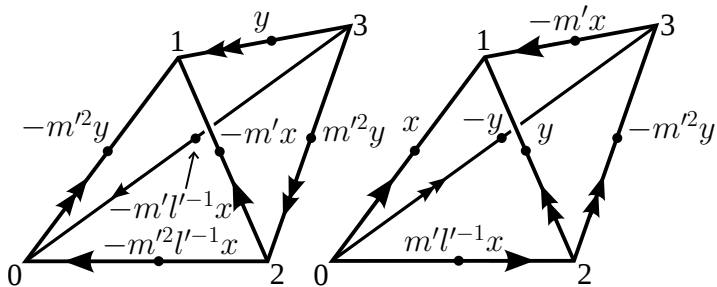
\begin{figure}[htb]
\scalebox{0.6}{\input{figures_gen/m004Sl2Extended.tex}}
\caption{Ptolemy coordinates, $n=2$.}\label{fig:PtolemyCoordinates_m004_n=2}
\end{figure}

The identification relations for the integral point corresponding to the single arrow are
\begin{equation}
\begin{aligned}
c_{12,0}&\overset{d}{=}c_{13,1},\quad &c_{13,1}&\overset{a}{=}(m^\prime)^{-1} l^\prime c_{02,0},\quad& c_{02,0}&\overset{c}{=}-m^\prime c_{02,1},\\
c_{02,1}&\overset{b}{=}-c_{03,0},\quad &c_{03,0}&\overset{a}{=}-m^\prime (l^\prime)^{-1} c_{01,1},\quad&c_{01,1}&\overset{c}{=}-(m^\prime)^{-1}c_{12,0}
\end{aligned}
\end{equation}
and for the double arrow they are
\begin{equation}
\begin{aligned}
c_{13,0}&\overset{b}{=}-c_{03,1},\quad &c_{03,1}&\overset{a}{=}-(m^\prime)^{-2}c_{23,0},\quad &c_{23,0}&\overset{d}{=}-c_{23,1}\\
c_{23,1}&\overset{b}{=}c_{01,0},\quad &c_{01,0}&\overset{c}{=}-(m^\prime)^2c_{12,1},\quad &c_{12,1}&\overset{d}{=}c_{13,0}.
\end{aligned}
\end{equation}
Here $\overset{\alpha}{=}$ indicates identification via $\alpha$. 
Letting $x=c_{01,1}$ and $y=c_{12,1}$, the Ptolemy coordinates are shown in Figure~\ref{fig:PtolemyCoordinates_m004_n=2}.

The Ptolemy relations
\begin{equation}
c_{03,k}c_{12,k}-c_{02,k}c_{13,k}+c_{01,k}c_{23,k},\qquad k=0,1
\end{equation}
become
\begin{equation}
m^{\prime2}(l^\prime)^{-1} x^2+m^{\prime 2}(l^\prime)^{-1} yx-m^{\prime 4}y^2=0,\qquad -y^2+m^{\prime 2}(l^\prime)^{-1} x^2-m^{\prime 2}xy
\end{equation}
or equivalently (using that $m^{\prime}=m$, $l^{\prime} =m^2l^{-1}$)
\begin{equation}
m^4lx^2+m^4lyx-m^4y^2=0,\qquad -y^2+m^4lx^2-m^2xy=0
\end{equation}
Plugging into magma, we obtain that the enhanced Ptolemy variety $\E P_n(\T)$ is given by
\begin{equation}
\begin{gathered}
x/y=\frac{m^8-m^6-m^4l-2m^4+1}{m^6 - m^2}\\
m^8l - m^6l - m^4 - 2m^4l - m^4l^2 - m^2l + l=0.
\end{gathered}
\end{equation}
We thus recover the known formula for the $\SL(2,\C)$ $A$-polynomial. The reduced Ptolemy variety $\E P_n(\T)$ is obtained by setting $y=1$.

\subsubsection{Recovering the representation}
One can explicitly compute the representation corresponding to a point in the Ptolemy variety using the face pairing presentation (see Section~\ref{sec:FacePairingPresentation}), but sometimes other presentations are more convenient.
The fundamental group of the figure $8$ knot complement also has a \emph{two-bridge} presentation  
\begin{equation}
\langle x_1,x_2\bigm\vert x_1w=wx_2, w=x_2x_1^{-1}x_2^{-1}x_1\rangle,
\end{equation}
and one easily checks that this presentation is isomorphic to the face pairing presentation~\ref{eq:FacePairingPresentation} via the isomorphism
\begin{equation}\label{eq:TwoBridgeFacePairing}
\begin{gathered}
a\mapsto x_2x_1x_2^{-1},\qquad b\mapsto x_1x_2^{-1},\qquad c\mapsto x_1\\
c\mapsto x_1,\qquad ab^{-1}\mapsfrom x_2.
\end{gathered}
\end{equation}
Using this, one obtains that the representation corresponding to a point in the enhanced Ptolemy variety is given by
\begin{equation}
x_1\mapsto \begin{pmatrix}m^{-1}&\frac{m^4l - m^2l - m^2 - l}{m^4 - 1}\\0&m\end{pmatrix},\qquad 
x_2\mapsto \begin{pmatrix}0&\frac{m^2l + m^2}{m^4 - 1}\\\frac{-m^8 + m^6 + 2m^4 + m^2l - 1}{m^6 - m^2}&\frac{m^2 + 1}{m}\end{pmatrix}.
\end{equation}
The explicit computations are unenlightening and have been omitted.

\subsubsection{The Dehn invariant}\label{sec:DehnInvariantExample}
Letting $z_1$ and $z_2$ denote the associated shapes, we have
\begin{equation}
z_1=\frac{x}{y},\qquad 1-z_1=-m^{\prime2}l^{\prime}\frac{y}{x},\qquad z_2=\frac{l^\prime}{m^{\prime2}}\frac{y^2}{x^2},\qquad 1-z_2=l^{\prime}\frac{y}{x}.
\end{equation}
The Dehn invariant $\alpha$ is thus given by 
\begin{equation}
\begin{aligned}
\alpha=&z_1\wedge(1-z_1)+z_2\wedge(1-z_2)=\frac{x}{y}\wedge\big(-m^{\prime2}l^{\prime}\frac{y}{x}\big)+\big(\frac{l^{\prime}}{m^{\prime2}}\frac{y^2}{x^2}\big)\wedge l^{\prime}\frac{y}{x}\\
=&-2m^\prime\wedge l^\prime\in\wedge^2(\C^*).
\end{aligned}
\end{equation}
Since $\overline m=m/m^{-1}=m^2$ and $\overline l=l/l^{-1}=l^2$ this agrees with Theorem~\ref{thm:DehnInvariantFormula}.

\subsection{Extended Ptolemy variety of the figure $8$ knot, $n=3$}
For $n=3$, the Ptolemy coordinates are shown in Figure~\ref{fig:PtolemyCoordinates_m004_n=3}.
The Ptolemy variety is thus given by
\begin{equation}
\begin{aligned}
-m_1m_2l_2x_1f_3-l_2m_2^2x_1f_1+m_1^2y_0f_2&=0,\qquad &y_0f_3m_1^3m_2l_1x_0f_2+x_0f_1&=0,\\
f_1m_2x_1+m_1f_3y_1-m_1m_2y_1f_0&=0, \qquad &f_2y_1-m_1m_2x_1f_3+f_0x_1&=0,\\
m_1l_1f_2x_0+m_2l_1x_0f_0-m_2y_1f_3&=0,\qquad &f_1y_0-m_1^3m_2^2l_1l_2x_1f_0+m_1^2m_2y_0f_3&=0,\\
-m_2m_1l_1x_0f_0-m_1^2l_1f_2y_0+m_2y_0f_1&=0,\qquad
&y_1f_0-f_1m_1x_0+m_1^2m_2y_1f_2&=0,
\end{aligned}
\end{equation}

\begin{figure}[htb]
\scalebox{0.9}{\input{figures_gen/ExtendedPtolemy_n=3.tex}}
\caption{Ptolemy coordinates, $n=3$.}\label{fig:PtolemyCoordinates_m004_n=3}
\end{figure}

The variety seems to be very difficult to compute. Using the fact that the geometric representation is the unique representation which is defined over $\Q(\sqrt{-3})$ and satisfies that $m_1\!=\!m_2\!=\!l_1\!=\!l_2\!=\!1$, one can use deformation techniques developed by Fabrice Rouillier (unpublished) to show that the geometric component is an algebraic covering with two leaves over the surface determined by
\begin{equation}\label{eq:Surface}
f_3^2m_1^2m_2-f_3m_1^2m_2^2-f_3m_1+f_3m_2+m_1m_2=0.
\end{equation}
The full description is rather large, so we shall not display it here. The computation was done by Fabrice Rouillier.

\section{Proof of Theorem~\ref{thm:EPnAndReps} and Theorem~\ref{thm:Projection}}
The proof of Theorems~\ref{thm:EPnAndReps} and~\ref{thm:Projection} follows the same strategy as the proof of the one-one correspondences~\eqref{eq:OneToOne} given in~\cite{GaroufalidisThurstonZickert}. The idea is to prove that a point in $\E P_n(\T)$ naturally determines a cocycle on $M$. The cocycle labels peripheral edges by elements in $B$ and thus determines a decorated (by $B$-cosets) representation. Conversely, a decorated representation determines a point in $\E P_n(\T)$. The diagonal action on $\E P_n(\T)$ corresponds to a coboundary action on cocycles, but does not change the decorated representation. The maps $p_i\colon \E P_n(\T)\to R_B(\partial M)$ are defined via the cocycle in such a way that Theorem~\ref{thm:Projection} is immediate from the construction.

\subsection{The face pairing presentation}\label{sec:FacePairingPresentation}
Given a fundamental polyhedron for $\T$ in $\widehat{\widetilde M}$, one has a presentation of $\pi_1(M)$ with a generator for each face pairing, and a relation for each $1$-cell (see Figures~\ref{fig:FacePairingGenerator} and \ref{fig:FacePairingRelation}). Note that the loop corresponding to a face pairing is the loop passing the face in the opposite direction.

\begin{figure}[htb]
\begin{minipage}[b]{0.48\textwidth}
\hspace{1cm}\scalebox{0.55}{\input{figures_gen/FacePairingGenerator.tex}}
\end{minipage}
\begin{minipage}[t]{0.48\textwidth}
\hspace{2cm}\scalebox{0.65}{\input{figures_gen/FacePairingRelation.tex}}
\end{minipage}
\\
\begin{minipage}[t]{0.48\textwidth}
\caption{Face pairing presentation, generator.}\label{fig:FacePairingGenerator}
\end{minipage}
\begin{minipage}[t]{0.48\textwidth}
\caption{Face pairing presentation, $edcba=1$.}\label{fig:FacePairingRelation}
\end{minipage}
\end{figure}

\begin{example}\label{ex:FacePairing}
For the triangulation of the figure 8 knot in Figure~\ref{fig:m004Triangulation}, one obtains a fundamental polyhedron, e.g.~by gluing together the faces corresponding to the face pairing $d$. By inspecting Figure~\ref{fig:m004Rectangle} one sees that the corresponding presentation is given by
\begin{equation}\label{eq:FacePairingPresentation}
\pi_1(M)=\big\langle a,b,c\bigm\vert ca^{-1}bc^{-1}a,ab^{-1}c^{-1}b\big\rangle.
\end{equation}
\end{example}
\subsection{Decorations}
We recall the basic properties of decorations and refer to \cite{GaroufalidisThurstonZickert} for more details. Let $\rho\colon\pi_1(M)\to\SL(n,\C)$ be boundary-Borel, and let $L=\widehat{\widetilde M}$. The triangulation $\T$ of $M$ induces a triangulation of $L$.
\begin{definition}
A \emph{decoration} of $\rho$ is a $\rho$-equivariant map $D\colon L^{(0)}\to \SL(n,\C)/B$, i.e.~a $\rho$-equivariant assignment of $B$-cosets to the vertices of $L$.
\end{definition}
Note that if $v$ is paired to $w$ via a face pairing $\alpha$, we must have $D(w)=\rho(\alpha)gB$  (see Figure~\ref{fig:FacePairingRelation}). The following is elementary.
\begin{lemma}\label{lemma:DecorationFreedom} Let $D$ be a decoration of $\rho$. If for some vertex $v$, $D(v)=gB$, we have
\begin{equation}\label{eq:DecorationFreedom}
g^{-1}\rho(\Stab(v))g\subset B,
\end{equation}
where $\Stab(v)\subset\pi_1(M)$ is the stabilizer of $v$.\qed 
\end{lemma}
\subsubsection{Freedom in the choice of decoration} The following is an elementary Corollary of Lemma~\ref{lemma:DecorationFreedom}. We shall not need it, so we leave its proof to the reader.
\begin{corollary} If $\rho$ maps each peripheral subgroup to a conjugate of $D$, and if each such contains an element with distinct eigenvalues, the set of decorations of $\rho$ is a torsor for $W^c$. Here $D$ is the group of diagonal matrices, $W$ is the Weyl group, and $c$ is the number of boundary components.\qed
\end{corollary}
\begin{remark}
Generally, there is more degree of freedom in the choice of decoration. If, for example, $\rho$ collapses the boundary component corresponding to $v$, then~\eqref{eq:DecorationFreedom} is satisfied for any $g$.
\end{remark}
\subsubsection{Peripheral holonomy}\label{sec:PeripheralHolonomy}
For each lift $v_i\in L^{(0)}$ of a boundary component $\partial_i M$ of $M$, we have a canonical isomorphism $\pi_1(\partial_i M)\cong\Stab(v_i)$.
\begin{definition}
The $i$th \emph{peripheral holonomy} of a decorated $\SL(n,\C)$-representation $\rho$ is the map
\begin{equation}
\xymatrix{\pi_1(\partial_i M)\cong\Stab(v_i)\ar[r]^-\rho&\rho(\Stab(v_i))\ar[r]&B}
\end{equation}
where the right map is conjugation by $D(v_i)$ (see~\eqref{eq:DecorationFreedom}), and $D$ is the decoration. It is defined up to conjugation by elements in $B$.
\end{definition}
Clearly, the peripheral holonomy only depends on $i$ and not on the choice of lift $v_i$. We thus have a \emph{restriction map}
\begin{equation}
r\colon\left\{\txt{Decorated, boundary-Borel\\$\pi_1(M)\to\SL(n,\C)$}\right\}\to R_B(\partial_1 M)\times\dots\times R_B(\partial_h M),
\end{equation}
where the $i$th component of $r$ takes a decorated representation to its $i$th peripheral holonomy.
\subsection{From Ptolemy coordinates to cocycles}
The construction of the cocycle hinges on the result below.
\begin{proposition}[Garoufalidis--Thurston--Zickert~\cite{GaroufalidisThurstonZickert}]\label{prop:NaturalCocycle}
\begin{equation}
\xymatrix{\left\{\txt{Ptolemy assign-\\ments on $\Delta^3_n$}\right\}\ar@{<->}[r]&\left\{\txt{Generic tuples\\$(g_0N,g_1N,g_2N,g_3N)$\\up to left action}\right\}\ar@{<->}[r]&\left\{\txt{Natural cocycles\\on $\overline\Delta^3_n$}\right\}}
\end{equation}
\end{proposition}
The Ptolemy assignment of a generic tuple $(g_0N,\dots,g_3N)$ is given by
\begin{equation}\label{eq:DecToPtolemy}
\Delta^3_n(\Z)\to\C^*, \qquad c_t=\det\big(\{g_0\}_{t_0}\cup\{g_1\}_{t_1}\cup\{g_2\}_{t_2}\cup\{g_3\}_{t_3}\big),
\end{equation}
and the natural cocycle labels long edges by counter diagonal elements $\alpha_{ij}$, and short edges by elements $\beta_{ijk}$ in $N$ (see Figure~\ref{fig:FattenedCocycle}). It is constructed using the fact that for (generic) $g_iN,g_jN$, there exist unique representatives $g_ix_i, g_jx_j$ such that $x_i^{-1}g_i^{-1}g_jx_j$ is counter diagonal. We refer to \cite{GaroufalidisThurstonZickert,GaroufalidisGoernerZickert} for explicit formulas in terms of the Ptolemy coordinates. The following is elementary.

\begin{lemma}\label{lemma:CocycleConjugation}
Let $\alpha_{ij}$ and $\beta_{ijk}$ be the natural cocycle for the tuple $(g_0N,\dots,g_3N)$. For diagonal matrices $d_0,\dots ,d_3$, let $\alpha_{ij}^\prime$ and $\beta_{ijk}^\prime$ be the cocycle of $(g_0d_0N,\dots,g_3d_3N)$. We then have
\begin{equation}
\pushQED{\qed} 
\alpha_{ij}^\prime=d_i^{-1}\alpha_{ij}d_j,\qquad \beta_{ijk}^\prime =d_i^{-1}\beta_{ijk}d_i^{-1}
\qedhere
\popQED
\end{equation}
\end{lemma}

An element in $\E P_n(\T)$ also determines a Ptolemy assignment on each simplex, but the corresponding cocycles no longer agree under the face pairings. 
To fix this, we shall consider a ``fattened'' decomposition of $M$, which mitigates this.
\subsubsection{The fattened decomposition}  Consider the polyhedral decomposition of $M$ obtained by thickening each hexagonal face, i.e.~replacing it with a hexagonal prism. An edge cycle ($1$-cell) of length $k$ is then replaced by a prism over a polygon with $k$ sides, so the decomposition consists of truncated simplices (one per simplex), hexagonal prisms (one per face), and polygonal prisms (one per edge). The top view of an edge link (as in Figure~\ref{fig:FacePairingRelation}) now looks like in Figure~\ref{fig:FattenedCuspTriangulation}. The triangulation has $3$ types of edges: \emph{long edges}, \emph{short edges}, and \emph{face-pairing-edges}.
\begin{definition}\label{def:NaturalCocycle}
The \emph{natural cocycle} associated to an element in $\E P_n(\T)$ is the cocycle labeling the long and short edges using Proposition~\ref{prop:NaturalCocycle}, and the face-pairing-edges as follows: If a face pairing $\alpha$ takes a vertex $v$ in $\T$ to $w$, the two face-pairing-edges near $v$ corresponding to $\alpha$ are both labeled by $M_{\alpha v}$ (see Figure~\ref{fig:FattenedCocycle}).
\end{definition}

\begin{figure}[htb]
\begin{minipage}[b]{0.48\textwidth}
\scalebox{0.55}{\input{figures_gen/FattenedCocycle.tex}}
\end{minipage}
\begin{minipage}[t]{0.38\textwidth}
\scalebox{0.7}{\input{figures_gen/FattenedCuspTriangulation.tex}}
\end{minipage}
\\
\begin{minipage}[t]{0.52\textwidth}
\caption{Hegagonal prism and corresponding cocycle.}\label{fig:FattenedCocycle}
\end{minipage}
\begin{minipage}[t]{0.46\textwidth}
\caption{Fattened version of Figure~\ref{fig:FacePairingRelation}.}\label{fig:FattenedCuspTriangulation}
\end{minipage}
\end{figure}

\begin{lemma}
The natural cocycle is indeed a cocycle, i.e.~the product of labelings around each face is $1$.
\end{lemma}
\begin{proof}
For the truncated simplices, this is immediate from Proposition~\ref{prop:NaturalCocycle}. We thus only need to check the following $3$ types of faces: polygons of face-pairing edges (as in the center of Figure~\ref{fig:FattenedCuspTriangulation}), squares consisting of two long edges and two face-pairing-edges, and squares consisting of two short edges and two face-pairing edges (see Figure~\ref{fig:FattenedCocycle}). For the polygons, this is obvious since each $M_\alpha$ is either $I$ or appears with its inverse. For the two types of squares, this follows from Lemma~\ref{lemma:CocycleConjugation}.
\end{proof}
\begin{corollary}\label{cor:PtolemyToRep}
A point in $\E P_n(\T)$ determines a decorated (by $B$-cosets) representation, which is invariant under the diagonal action on $\E P_n(\T)$.\qed
\end{corollary}
\begin{proof}
By construction, the diagonal action on $\E P_n(\T)$ corresponds to the action by coboundaries on natural cocycles. This neither changes the representation, nor the $B$-cosets of the decoration.
\end{proof}

\subsection{From decorated representations to Ptolemy coordinates}
Let $\rho$ be a decorated boundary-Borel representation. Pick fundamental rectangles $R_i$, and let $\mu_i^\prime$ and $\lambda_i^\prime$ be the corresponding peripheral curves (see Section~\ref{sec:ExtendedPtolemy}). Define diagonal matrices $D_{\mu_i^\prime}=\diag(m_{i,1}^\prime,\dots,m_{h,1}^\prime)$ and $D_{\lambda_i^\prime}=\diag(l_{i,1}^\prime,\dots,l_{h,1}^\prime)$ to be the diagonal part of the $i$th boundary homology of $\mu_i^\prime$ and $\lambda_i^{\prime}$, respectively.
The decoration assigns to each triangle $t$ in each $R_i$ a $B$-coset $D(t)$. The following is clear from the definitions.
\begin{lemma} There exists a lift
\begin{equation}
\widetilde D\colon L^{(0)}\to\SL(n,\C)/N,
\end{equation}
of $D$ satisfying that when triangles $t$ and $t^\prime$ are paired via a face pairing $\alpha$, we have
\begin{equation}\label{eq:LiftToNCosets}
\rho(\alpha)\widetilde D(t)M_{\alpha v_t}=\widetilde D(t^\prime)
\end{equation}
where $v_t$ is the vertex index corresponding to $t$, and $M_{\alpha}$ is defined as in Section~\ref{sec:ExtendedPtolemy}. Moreover, any two such lifts differ by the action of $H^c$ by right multiplication.\qed
\end{lemma}
We thus obtain a decoration by $N$-cosets, hence a Ptolemy assignment, on each simplex, and \eqref{eq:LiftToNCosets} implies that the Ptolemy coordinates satisfy the identification relations. We thus have.
\begin{corollary}\label{cor:RepToPtolemy}
A decorated boundary-Borel representation determines an element in $\E P_n(\T)_{\red}$.\qed
\end{corollary}

It is clear from the constructions that the maps in Corollaries~\ref{cor:PtolemyToRep} and \ref{cor:RepToPtolemy} are inverses. Also note that if we use a different choice $\{R_i^\prime\}$ of fundamental rectangles, the map 
\begin{equation}
\E P_n(\T)^{\{R_i\}}_{\red}\to \E P_n(\T)^{\{R_i^\prime\}}_{\red}
\end{equation}
induced by Corollaries~\ref{cor:PtolemyToRep} and \ref{cor:RepToPtolemy} is a regular bijection. Hence, the choice of fundamental rectangles is inessential.

\section{Proof of Theorem~\ref{thm:DehnInvariantFormula}}
Theorem~\ref{thm:DehnInvariantFormula} is a corollary of the symplectic properties of the generalized gluing equations proved in~\cite{GaroufalidisZickert} (see also~\cite{Guilloux}) and elementary properties of shape assignments developed in~\cite{GaroufalidisGoernerZickert}.
\subsection{Shape assignments}\label{sec:ShapeAssignments}
\begin{definition}
A \emph{subsimplex} of $\Delta^3_n$ is a translate of $\Delta^3_2\subset\R^4$ by an element $s\in\Delta_{n-2}(\Z)$. 
\end{definition}
\begin{definition}
A \emph{shape assignment} on $\Delta_n^3$ is a map
\begin{equation}
\Delta_n^3(\Z)\times\dot\Delta^3_2(\Z)\to\C\setminus\{0,1\}, \qquad (s,e)\mapsto z_s^e
\end{equation}
satisfying the relations
\begin{equation}\label{eq:ShapeRelations}
z_s^{1100}=z_s^{0011}=\frac{1}{1-z_s^{0101}},\quad z_s^{0110}=z_s^{1001}=\frac{1}{1-z_s^{0011}},\quad z_s^{1010}=z_s^{0101}=\frac{1}{1-z_s^{0110}}.
\end{equation}
A \emph{shape assignment} on $\T$ is a shape assignment $\{z_{s,\Delta}^e\}$ on each simplex such that for each integral point $p=\{(t,\Delta)\mid t\in\Delta_n(\Z),\Delta\in\T\}$, the \emph{shape parameters} $z_{s,\Delta}^e$ satisfy the \emph{generalized gluing equation}
\begin{equation}\label{eq:GeneralizedGluingEquations}
\prod_{(t,\Delta)\in p}\,\prod_{t=s+e}z^e_{s,\Delta}=1.
\end{equation}
The variety of shape assignments on $\T$ is denoted by $V_n(\T)$. 
\end{definition}

\begin{remark} The set $\dot\Delta^3_2(\Z)$ parametrizes the (unoriented) edges of an ordered simplex, so one can think of a shape assignment on $\Delta^3_n$ as an assignment of a shape parameter to each edge of each \emph{subsimplex} of $\Delta^3_n$.
\end{remark}

\subsection{The generalized Neumann chain complex}
For an ordered simplex $\Delta$, let $J_\Delta$ be the free abelian group on the (unoriented) edges $\varepsilon_{ij}$ of $\Delta$ subject to the relations
\begin{equation}
\varepsilon_{01}=\varepsilon_{23},\qquad \varepsilon_{12}=\varepsilon_{03},\qquad \varepsilon_{02}=\varepsilon_{13},\qquad \varepsilon_{01}+\varepsilon_{12}+\varepsilon_{02}=0.
\end{equation}
We endow $J_\Delta$ with the non-degenerate skew symmetric bilinear form $\Omega$ defined by
\begin{equation}
\Omega(\varepsilon_{01},\varepsilon_{12})=\Omega(\varepsilon_{12},\varepsilon_{02})=\Omega(\varepsilon_{02},\varepsilon_{01})=1.
\end{equation}
Note that in the basis $\{\varepsilon_{01},\varepsilon_{12}\}$ for $J_\Delta$, $\Omega$ is the standard symplectic form given by $\left(\begin{smallmatrix}0&1\\-1&0\end{smallmatrix}\right)$.
In~\cite{GaroufalidisZickert} we defined a chain complex $\J^n$ given by
\begin{equation}
\xymatrix{0\ar[r]&C_0^n(\T)\ar[r]^-\alpha&C_1^n(\T)\ar[r]^-\beta&J^n(\T)\ar[r]^{\beta^*}&C_1^n(\T)\ar[r]^-{\alpha^*}&C_0^n(\T)\ar[r]&0}.
\end{equation}
Here $C_1^n(\T)$ is the free abelian group on integral points and $J^n(\T)=\oplus_{\Delta\in\T} J_\Delta$ with $\Omega$ extended orthogonally. The map $\beta$ encodes the gluing equations and $\beta^*$ encodes the Ptolemy relations. We shall not need the precise definitions. 
The form $\Omega$ descends to a form on $H_3(\J^n)$, which is non-degenerate modulo torsion.
\begin{theorem}[Garoufalidis--Zickert~\cite{GaroufalidisZickert}]\label{thm:SympProperty}
We have an isomorphism
\begin{equation}
\delta^{\prime}\colon H_1(\partial M;\Z[1/n]^{n-1})\cong H_3(\J^n)\otimes\Z[1/n],
\end{equation}
such that the form $\Omega$ corresponds to the form $\omega_A$ on $H_1(\partial M;\Z[1/n]^{n-1})$ given by
\begin{equation}
\omega_A(\alpha\otimes v,\beta\otimes v)=\iota(\alpha,\beta)\langle v,Aw\rangle.
\end{equation}
Here, $A$ is the Cartan matrix for $\SL(n,\C)$, $\iota$ is the intersection form on $H_1(\partial M)$, and $\langle-,-\rangle$ is the standard inner product on $\Z^{n-1}$.\qed
\end{theorem}
Let $\mu\subset \C^*$ be the roots of unity. The following is an elementary corollary of Theorem~\ref{thm:SympProperty}.
\begin{corollary}
Under the dual isomorphism 
\begin{equation}\label{eq:DeltaPrimeIso}
\delta^{\prime *}\colon H^3(\J^n)\otimes\C^*/\mu\cong H^1(\partial M;(\C^*/\mu)^{n-1}),
\end{equation}
the form $\Omega^*$ corresponds to the form $\omega_A^*$ given by
\begin{equation}\label{eq:OmegaStarCartan}
\pushQED{\qed}
\omega_A^*(\alpha\otimes x,\beta\otimes y)=\alpha\cup\beta([\partial M])\langle x,A^{-1}y\rangle.\qedhere
\popQED
\end{equation}
\end{corollary}
\subsection{The proof}
The idea is that a shape assignment $z\in V_n(\T)$ determines two elements 
\begin{equation}
J(z)\in H^3(\J;\C^*/\mu),\qquad C(z)\in H^1(\partial M;\C^*/\mu),
\end{equation} which agree under the isomorphism~\eqref{eq:DeltaPrimeIso}.
The result is then an elementary computation using~\eqref{eq:OmegaStarCartan} and some basic properties of bilinear forms and exterior powers.
\subsubsection{Definition of $J(z)$} A shape assignment assigns a shape parameter in $\C^*$ to each edge of each subsimplex, and can thus be viewed as an element $J(z)$ in $\Hom(J^n(\T),\C^*/\pm 1)$; the sign is due to the fact that $z_s^{1100}z_s^{0110}z_s^{1010}=-1$ (see~\eqref{eq:ShapeRelations}) and not $1$. Since the shape parameters satisfy the gluing equations, $J(z)$ a cocycle, so determines a cohomology class in $H^3(\J,\C^*/\pm 1)$.
\subsubsection{Definition of $C(z)$}  By the main result of~\cite{GaroufalidisGoernerZickert}, a shape assignment $z$ determines a generically decorated $\PGL(n,\C)$-representation (and vice versa), which is given explicitly by a cocycle on $M$. The cocycle labels peripheral edges by elements in $B\subset\PGL(n,\C)$, so every edge path in $\partial M$ has a well defined holonomy in $B$ (the peripheral holonomy; Section~\ref{sec:PeripheralHolonomy}). Composing with the map $B\to(\C^*)^{n-1}$ taking ratios of consecutive diagonal entries, we obtain a cohomology class
\begin{equation}
C(z)\in H^1(\partial M,(\C^*)^{n-1}).
\end{equation}
By construction, this class is trivial if and only if the corresponding representation is boundary-unipotent.

\begin{theorem}[Garoufalidis--Zickert~\cite{GaroufalidisZickert}]\label{thm:JandCAgree}
The elements $J(z)$ and $C(z)$ agree under the the isomorphism~\eqref{eq:DeltaPrimeIso}, i.e.~$\delta^{\prime*}(J(z))=C(z)$ (in fact, this holds in $H^*(\partial M;(\C^*/\pm 1)^{n-1})$ as well).\qed
\end{theorem}

\subsubsection{Bilinear forms and exterior powers} Let $A$ be an abelian group with a bilinear form $\omega\colon A\times A\to\Z$. As noted in~\cite{BergeronFalbelGuilloux}, we have for each abelian group $B$ (we shall only consider $B=\C^*/\mu$) a bilinear map
\begin{equation}
\wedge_\omega\colon (A\otimes_\Z B)\times (A\otimes_\Z B)\to \wedge^2_\Z(B), \qquad (a\otimes v,b\otimes w)\mapsto \omega(a,b)v\wedge w.
\end{equation}
\begin{lemma} We have
\begin{equation}\label{eq:FundamentalEq}
J(z)\wedge_{\Omega^*}J(z)=C(z)\wedge_{\omega_A^*} C(z)\in\wedge^2_\Z(\C^*/\mu).
\end{equation}
\end{lemma}
\begin{proof}
This follows from Theorem~\ref{thm:JandCAgree} and~\eqref{eq:OmegaStarCartan} using basic functoriality properties.  
\end{proof}

\subsubsection{Basic computations}
We now compute the left and right hand side of~\eqref{eq:FundamentalEq}.
Letting $\mu_i^*$ and $\lambda_i^*$ in $H^1(\partial M)$ denote the duals of $\mu_i$ and $\lambda_i$, respectively, we can write $C(z)$ in the form 
\begin{equation}\label{eq:ConcreteCformula}
C(z)=\sum_i(\mu_i\otimes(\overline m_{i,1}(z),\dots,\overline m_{i,n-1}(z))+\lambda_i\otimes(\overline l_{i,1}(z),\dots,\overline l_{i,n-1}(z)),
\end{equation}
where $\overline m_{i,j}(z)=m_{i,j}(z)/m_{i,j+1}(z)$, and $m_{i,j}(z)$ is the $j$th diagonal entry in the boundary holonomy of $\mu_i$ ($\overline l_{i,j}(z)$ is defined similarly). 
\begin{lemma} We have
\begin{equation}
C(z)\wedge_{\omega_A^*}C(z)=2\sum_{i=1}^h\sum_{j,k}A^{-1}_{jk}\overline m_{i,j}(z)\wedge \overline l_{i,k}(z)\in\wedge^2_\Z(\C^*/\mu)
\end{equation}
\end{lemma}
\begin{proof}
Using~\eqref{eq:ConcreteCformula} this follows immediately from the formula~\eqref{eq:OmegaStarCartan}.
\end{proof}

The following is shown for $n=3$ in~\cite{BergeronFalbelGuilloux} (the missing minus sign there is due to a different cross-ratio convention).
\begin{lemma}
The element $J(z)\wedge_{\Omega^*}J(z)$ is minus $2$ times the Dehn invariant of $z$.
\end{lemma}
\begin{proof}
Since $J$ is a direct sum of $J_\Delta$'s, it is enough to prove this for a shape assignment on $\Delta_2^3$. We denote $z^{1100}$ and $z^{0110}$ by $z$ and $z^\prime$, respectively ($z$ denotes both the shape assignment and the shape).
Using the basis $\{\varepsilon_{01},\varepsilon_{12}\}$ for $J_\Delta$, and the corresponding dual basis $\{\varepsilon_{01}^*,\varepsilon_{12}^*\}$ for $J_\Delta^*$, one has 
\begin{equation}
J(z)=\varepsilon_{01}^*\otimes z+\varepsilon_{12}^*\otimes z^\prime.
\end{equation}
Since the $\Omega$-duals of $\varepsilon_{01}$ and $\varepsilon_{12}$ are $-\varepsilon_{12}$ and $\varepsilon_{01}$, respectively, we have
\begin{equation}
J(z)\wedge_{\Omega^*}J(z)=\Omega(\varepsilon_{01}\otimes z^\prime-\varepsilon_{12}\otimes z,\varepsilon_{01}\otimes z^\prime-\varepsilon_{12}\otimes z)
=-z^{\prime}\wedge z+z\wedge z^\prime=-2z\wedge (1-z), 
\end{equation}
where the last equation follows from~\eqref{eq:ShapeRelations}. This proves the result.
\end{proof}

Since $\wedge^2_\Z(\C^*)=\wedge^2_\Z(\C^*/\mu)$, we have proved.
\begin{theorem}
The Dehn invariant of a shape assignment $z\in V_n(\T)$ is given by
\begin{equation}
\pushQED{\qed}
\nu\circ\lambda(z)=-\sum_{i=1}^h\sum_{j,k}A^{-1}_{jk}\overline m_{i,j}(z)\wedge \overline l_{i,k}(z)\in\wedge^2_\Z(\C^*).\qedhere
\popQED
\end{equation}
\end{theorem}

Theorem~\ref{thm:DehnInvariantFormula} is an immediate consequence.

\begin{remark}
When working over other fields one needs to be more careful since $\wedge^2(k^*)$ is not $n$-divisible in general. There are also torsion issues due to vertex orderings. We shall not discuss these issues here.
\end{remark}

\bibliographystyle{plain}
\bibliography{BibFile}

\end{document}

%% file: figures_gen/IdentificationRelations_n=3.tex
\begingroup
 \setlength{\unitlength}{0.8pt}
 \begin{picture}(254.00658,129.01921)
 \put(0,0){\includegraphics{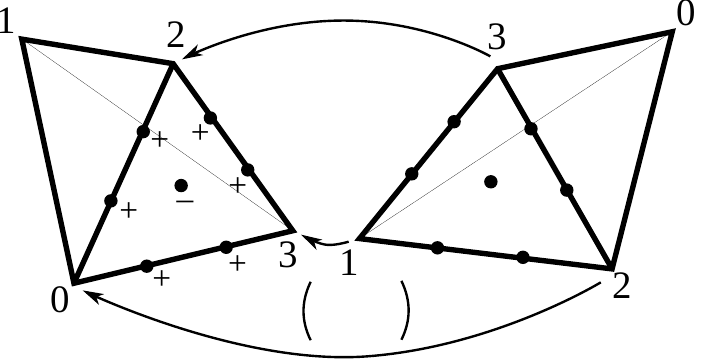}}

\definecolor{inkcol1}{rgb}{0.0,0.0,0.0}
   \put(20.800672,91.337226){\rotatebox{360.0}{\makebox(0,0)[tl]{\strut{}{
    \begin{minipage}[h]{127.655896pt}
\textcolor{inkcol1}{\large{$\Delta_0$}}\\
\end{minipage}}}}}%

\definecolor{inkcol1}{rgb}{0.0,0.0,0.0}
   \put(203.561002,93.785866){\rotatebox{360.0}{\makebox(0,0)[tl]{\strut{}{
    \begin{minipage}[h]{127.655896pt}
\textcolor{inkcol1}{\large{$\Delta_1$}}\\
\end{minipage}}}}}%

\definecolor{inkcol1}{rgb}{0.0,0.0,0.0}
   \put(109.959922,27.701426){\rotatebox{360.0}{\makebox(0,0)[tl]{\strut{}{
    \begin{minipage}[h]{127.655896pt}
\textcolor{inkcol1}{\normalsize{$0123$}}\\
\end{minipage}}}}}%

\definecolor{inkcol1}{rgb}{0.0,0.0,0.0}
   \put(109.959922,17.094826){\rotatebox{360.0}{\makebox(0,0)[tl]{\strut{}{
    \begin{minipage}[h]{127.655896pt}
\textcolor{inkcol1}{\normalsize{$1302$}}\\
\end{minipage}}}}}%

\definecolor{inkcol1}{rgb}{0.0,0.0,0.0}
   \put(78.446002,23.604616){\rotatebox{360.0}{\makebox(0,0)[tl]{\strut{}{
    \begin{minipage}[h]{127.655896pt}
\textcolor{inkcol1}{\normalsize{$\sigma\!=$}}\\
\end{minipage}}}}}%

 \end{picture}
\endgroup

%% file: figures_gen/IdentificationRelations_n=4.tex
\begingroup
 \setlength{\unitlength}{0.8pt}
 \begin{picture}(254.00659,129.01921)
 \put(0,0){\includegraphics{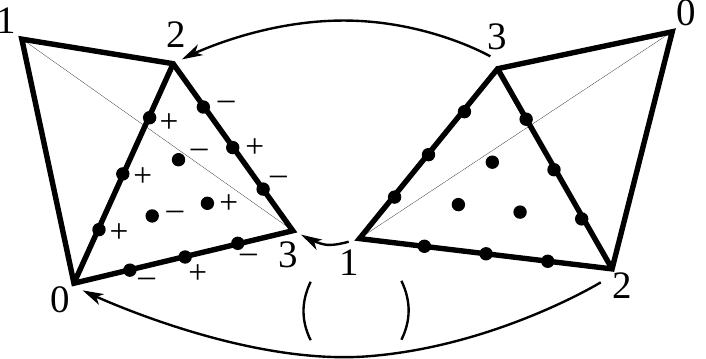}}

\definecolor{inkcol1}{rgb}{0.0,0.0,0.0}
   \put(20.8006767,91.337226){\rotatebox{360.0}{\makebox(0,0)[tl]{\strut{}{
    \begin{minipage}[h]{127.655896pt}
\textcolor{inkcol1}{\large{$\Delta_0$}}\\
\end{minipage}}}}}%

\definecolor{inkcol1}{rgb}{0.0,0.0,0.0}
   \put(203.5610098,93.785866){\rotatebox{360.0}{\makebox(0,0)[tl]{\strut{}{
    \begin{minipage}[h]{127.655896pt}
\textcolor{inkcol1}{\large{$\Delta_1$}}\\
\end{minipage}}}}}%

\definecolor{inkcol1}{rgb}{0.0,0.0,0.0}
   \put(66.5723698,107.415642){\rotatebox{360.0}{\makebox(0,0)[tl]{\strut{}{
    \begin{minipage}[h]{93.708112pt}
\textcolor{inkcol1}{\normalsize{$t_0\!\!=\!0031$}}\\
\end{minipage}}}}}%

\definecolor{inkcol1}{rgb}{0.0,0.0,0.0}
   \put(104.5693701,94.051042){\rotatebox{360.0}{\makebox(0,0)[tl]{\strut{}{
    \begin{minipage}[h]{90.982816pt}
\textcolor{inkcol1}{\normalsize{$t_1\!\!=\!0103$}}\\
\end{minipage}}}}}%

\definecolor{inkcol1}{rgb}{0.0,0.0,0.0}
   \put(109.9599268,27.701426){\rotatebox{360.0}{\makebox(0,0)[tl]{\strut{}{
    \begin{minipage}[h]{127.655896pt}
\textcolor{inkcol1}{\normalsize{$0123$}}\\
\end{minipage}}}}}%

\definecolor{inkcol1}{rgb}{0.0,0.0,0.0}
   \put(109.9599268,17.094826){\rotatebox{360.0}{\makebox(0,0)[tl]{\strut{}{
    \begin{minipage}[h]{127.655896pt}
\textcolor{inkcol1}{\normalsize{$1302$}}\\
\end{minipage}}}}}%

\definecolor{inkcol1}{rgb}{0.0,0.0,0.0}
   \put(78.4460088,23.604616){\rotatebox{360.0}{\makebox(0,0)[tl]{\strut{}{
    \begin{minipage}[h]{127.655896pt}
\textcolor{inkcol1}{\normalsize{$\sigma\!=$}}\\
\end{minipage}}}}}%

 \end{picture}
\endgroup

%% file: figures_gen/DiagonalAction.tex
\begingroup
 \setlength{\unitlength}{0.8pt}
 \begin{picture}(645.39752,196.80019)
 \put(0,0){\includegraphics{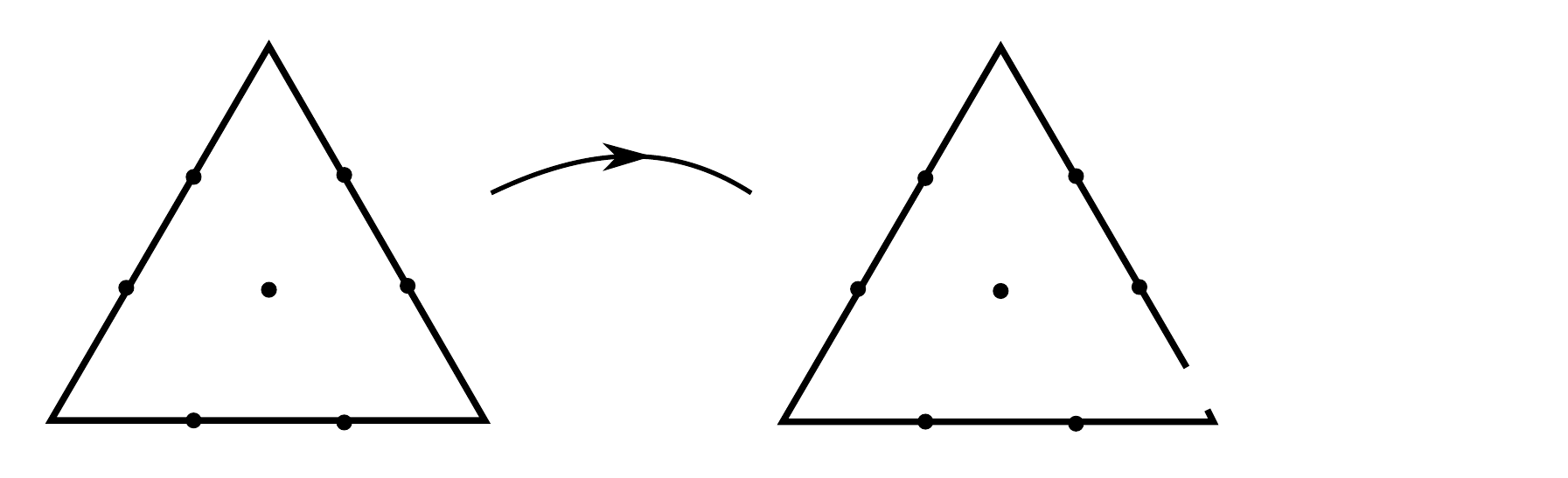}}

\definecolor{inkcol1}{rgb}{0.0,0.0,0.0}
   \put(62.641013,144.164248){\rotatebox{360.0}{\makebox(0,0)[tl]{\strut{}{
    \begin{minipage}[h]{129.560208pt}
\textcolor{inkcol1}{\LARGE{$a$}}\\
\end{minipage}}}}}%

\definecolor{inkcol1}{rgb}{0.0,0.0,0.0}
   \put(99.460543,204.189388){\rotatebox{360.0}{\makebox(0,0)[tl]{\strut{}{
    \begin{minipage}[h]{129.560208pt}
\textcolor{inkcol1}{\LARGE{$g_1N$}}\\
\end{minipage}}}}}%

\definecolor{inkcol1}{rgb}{0.0,0.0,0.0}
   \put(-0.738019,18.101318){\rotatebox{360.0}{\makebox(0,0)[tl]{\strut{}{
    \begin{minipage}[h]{129.560208pt}
\textcolor{inkcol1}{\LARGE{$g_0N$}}\\
\end{minipage}}}}}%

\definecolor{inkcol1}{rgb}{0.0,0.0,0.0}
   \put(141.150743,148.560638){\rotatebox{360.0}{\makebox(0,0)[tl]{\strut{}{
    \begin{minipage}[h]{129.560208pt}
\textcolor{inkcol1}{\LARGE{$b$}}\\
\end{minipage}}}}}%

\definecolor{inkcol1}{rgb}{0.0,0.0,0.0}
   \put(36.094863,100.073308){\rotatebox{360.0}{\makebox(0,0)[tl]{\strut{}{
    \begin{minipage}[h]{129.560208pt}
\textcolor{inkcol1}{\LARGE{$c$}}\\
\end{minipage}}}}}%

\definecolor{inkcol1}{rgb}{0.0,0.0,0.0}
   \put(107.310623,104.113918){\rotatebox{360.0}{\makebox(0,0)[tl]{\strut{}{
    \begin{minipage}[h]{129.560208pt}
\textcolor{inkcol1}{\LARGE{$d$}}\\
\end{minipage}}}}}%

\definecolor{inkcol1}{rgb}{0.0,0.0,0.0}
   \put(169.434993,100.073308){\rotatebox{360.0}{\makebox(0,0)[tl]{\strut{}{
    \begin{minipage}[h]{129.560208pt}
\textcolor{inkcol1}{\LARGE{$e$}}\\
\end{minipage}}}}}%

\definecolor{inkcol1}{rgb}{0.0,0.0,0.0}
   \put(65.894363,51.080918){\rotatebox{360.0}{\makebox(0,0)[tl]{\strut{}{
    \begin{minipage}[h]{129.560208pt}
\textcolor{inkcol1}{\LARGE{$f$}}\\
\end{minipage}}}}}%

\definecolor{inkcol1}{rgb}{0.0,0.0,0.0}
   \put(139.635503,47.040308){\rotatebox{360.0}{\makebox(0,0)[tl]{\strut{}{
    \begin{minipage}[h]{129.560208pt}
\textcolor{inkcol1}{\LARGE{$g$}}\\
\end{minipage}}}}}%

\definecolor{inkcol1}{rgb}{0.0,0.0,0.0}
   \put(195.851243,23.502498){\rotatebox{360.0}{\makebox(0,0)[tl]{\strut{}{
    \begin{minipage}[h]{129.560208pt}
\textcolor{inkcol1}{\LARGE{$g_2N$}}\\
\end{minipage}}}}}%

\definecolor{inkcol1}{rgb}{0.0,0.0,0.0}
   \put(628.742383,165.496518){\rotatebox{360.0}{\makebox(0,0)[tl]{\strut{}{
    \begin{minipage}[h]{129.560208pt}
\textcolor{inkcol1}{\LARGE{$\diag(x_1,x_2,x_3)$}}\\
\end{minipage}}}}}%

\definecolor{inkcol1}{rgb}{0.0,0.0,0.0}
   \put(627.612053,117.419378){\rotatebox{360.0}{\makebox(0,0)[tl]{\strut{}{
    \begin{minipage}[h]{129.560208pt}
\textcolor{inkcol1}{\LARGE{$\diag(y_1,y_2,y_3)$}}\\
\end{minipage}}}}}%

\definecolor{inkcol1}{rgb}{0.0,0.0,0.0}
   \put(628.622213,68.426988){\rotatebox{360.0}{\makebox(0,0)[tl]{\strut{}{
    \begin{minipage}[h]{129.560208pt}
\textcolor{inkcol1}{\LARGE{$\diag(z_1,z_2,z_3)$}}\\
\end{minipage}}}}}%

\definecolor{inkcol1}{rgb}{0.0,0.0,0.0}
   \put(565.101353,163.051928){\rotatebox{360.0}{\makebox(0,0)[tl]{\strut{}{
    \begin{minipage}[h]{129.560208pt}
\textcolor{inkcol1}{\LARGE{$D_0=$}}\\
\end{minipage}}}}}%

\definecolor{inkcol1}{rgb}{0.0,0.0,0.0}
   \put(565.055673,115.669658){\rotatebox{360.0}{\makebox(0,0)[tl]{\strut{}{
    \begin{minipage}[h]{129.560208pt}
\textcolor{inkcol1}{\LARGE{$D_1=$}}\\
\end{minipage}}}}}%

\definecolor{inkcol1}{rgb}{0.0,0.0,0.0}
   \put(566.065833,66.172168){\rotatebox{360.0}{\makebox(0,0)[tl]{\strut{}{
    \begin{minipage}[h]{129.560208pt}
\textcolor{inkcol1}{\LARGE{$D_2=$}}\\
\end{minipage}}}}}%

\definecolor{inkcol1}{rgb}{0.0,0.0,0.0}
   \put(303.172203,149.725778){\rotatebox{360.0}{\makebox(0,0)[tl]{\strut{}{
    \begin{minipage}[h]{129.560208pt}
\textcolor{inkcol1}{\LARGE{$ax_1y_1y_2$}}\\
\end{minipage}}}}}%

\definecolor{inkcol1}{rgb}{0.0,0.0,0.0}
   \put(290.138833,18.398978){\rotatebox{360.0}{\makebox(0,0)[tl]{\strut{}{
    \begin{minipage}[h]{129.560208pt}
\textcolor{inkcol1}{\LARGE{$g_0D_0N$}}\\
\end{minipage}}}}}%

\definecolor{inkcol1}{rgb}{0.0,0.0,0.0}
   \put(384.220073,203.256908){\rotatebox{360.0}{\makebox(0,0)[tl]{\strut{}{
    \begin{minipage}[h]{129.560208pt}
\textcolor{inkcol1}{\LARGE{$g_1D_1N$}}\\
\end{minipage}}}}}%

\definecolor{inkcol1}{rgb}{0.0,0.0,0.0}
   \put(489.781013,25.470048){\rotatebox{360.0}{\makebox(0,0)[tl]{\strut{}{
    \begin{minipage}[h]{129.560208pt}
\textcolor{inkcol1}{\LARGE{$g_2D_2N$}}\\
\end{minipage}}}}}%

\definecolor{inkcol1}{rgb}{0.0,0.0,0.0}
   \put(275.448483,101.267718){\rotatebox{360.0}{\makebox(0,0)[tl]{\strut{}{
    \begin{minipage}[h]{129.560208pt}
\textcolor{inkcol1}{\LARGE{$cx_1x_2y_1$}}\\
\end{minipage}}}}}%

\definecolor{inkcol1}{rgb}{0.0,0.0,0.0}
   \put(444.769023,147.987278){\rotatebox{360.0}{\makebox(0,0)[tl]{\strut{}{
    \begin{minipage}[h]{129.560208pt}
\textcolor{inkcol1}{\LARGE{$by_1y_2z_1$}}\\
\end{minipage}}}}}%

\definecolor{inkcol1}{rgb}{0.0,0.0,0.0}
   \put(371.075573,104.099478){\rotatebox{360.0}{\makebox(0,0)[tl]{\strut{}{
    \begin{minipage}[h]{129.560208pt}
\textcolor{inkcol1}{\LARGE{$dx_1y_1z_1$}}\\
\end{minipage}}}}}%

\definecolor{inkcol1}{rgb}{0.0,0.0,0.0}
   \put(472.732053,96.520258){\rotatebox{360.0}{\makebox(0,0)[tl]{\strut{}{
    \begin{minipage}[h]{129.560208pt}
\textcolor{inkcol1}{\LARGE{$ey_1z_1z_2$}}\\
\end{minipage}}}}}%

\definecolor{inkcol1}{rgb}{0.0,0.0,0.0}
   \put(335.019543,50.352188){\rotatebox{360.0}{\makebox(0,0)[tl]{\strut{}{
    \begin{minipage}[h]{129.560208pt}
\textcolor{inkcol1}{\LARGE{$fx_1x_2z_1$}}\\
\end{minipage}}}}}%

\definecolor{inkcol1}{rgb}{0.0,0.0,0.0}
   \put(426.026083,44.176178){\rotatebox{360.0}{\makebox(0,0)[tl]{\strut{}{
    \begin{minipage}[h]{129.560208pt}
\textcolor{inkcol1}{\LARGE{$gx_1z_1z_2$}}\\
\end{minipage}}}}}%

 \end{picture}
\endgroup

%% file: figures_gen/FacePairing.tex
\begingroup
 \setlength{\unitlength}{0.8pt}
 \begin{picture}(460.80347,252.21315)
 \put(0,0){\includegraphics{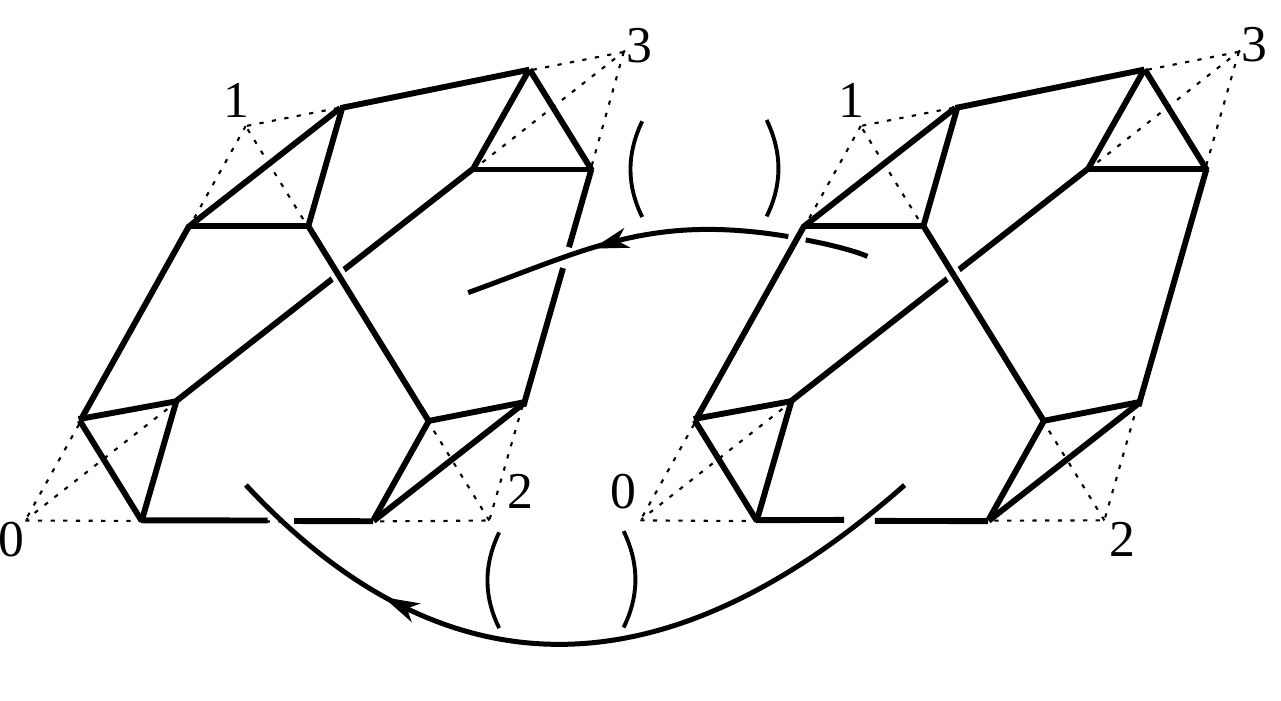}}

\definecolor{inkcol1}{rgb}{0.0,0.0,0.0}
   \put(195.39645,17.311006){\rotatebox{360.0}{\makebox(0,0)[tl]{\strut{}{
    \begin{minipage}[h]{127.655896pt}
\textcolor{inkcol1}{\Large{$a$}}\\
\end{minipage}}}}}%

\definecolor{inkcol1}{rgb}{0.0,0.0,0.0}
   \put(230.7518,205.303996){\rotatebox{360.0}{\makebox(0,0)[tl]{\strut{}{
    \begin{minipage}[h]{127.655896pt}
\textcolor{inkcol1}{\Large{$0123$}}\\
\end{minipage}}}}}%

\definecolor{inkcol1}{rgb}{0.0,0.0,0.0}
   \put(230.7518,189.646626){\rotatebox{360.0}{\makebox(0,0)[tl]{\strut{}{
    \begin{minipage}[h]{127.655896pt}
\textcolor{inkcol1}{\Large{$2103$}}\\
\end{minipage}}}}}%

\definecolor{inkcol1}{rgb}{0.0,0.0,0.0}
   \put(179.26157,57.319186){\rotatebox{360.0}{\makebox(0,0)[tl]{\strut{}{
    \begin{minipage}[h]{127.655896pt}
\textcolor{inkcol1}{\Large{$0123$}}\\
\end{minipage}}}}}%

\definecolor{inkcol1}{rgb}{0.0,0.0,0.0}
   \put(179.26157,41.661816){\rotatebox{360.0}{\makebox(0,0)[tl]{\strut{}{
    \begin{minipage}[h]{127.655896pt}
\textcolor{inkcol1}{\Large{$0213$}}\\
\end{minipage}}}}}%

\definecolor{inkcol1}{rgb}{0.0,0.0,0.0}
   \put(243.72311,165.720086){\rotatebox{360.0}{\makebox(0,0)[tl]{\strut{}{
    \begin{minipage}[h]{127.655896pt}
\textcolor{inkcol1}{\Large{$b$}}\\
\end{minipage}}}}}%

\definecolor{inkcol1}{rgb}{0.0,0.0,0.0}
   \put(83.69421,51.141356){\rotatebox{360.0}{\makebox(0,0)[tl]{\strut{}{
    \begin{minipage}[h]{127.655896pt}
\textcolor{inkcol1}{\LARGE{$\Delta_i$}}\\
\end{minipage}}}}}%

\definecolor{inkcol1}{rgb}{0.0,0.0,0.0}
   \put(316.13646,54.357786){\rotatebox{360.0}{\makebox(0,0)[tl]{\strut{}{
    \begin{minipage}[h]{127.655896pt}
\textcolor{inkcol1}{\LARGE{$\Delta_j$}}\\
\end{minipage}}}}}%

 \end{picture}
\endgroup

%% file: figures_gen/CuspImages.tex
\begingroup
 \setlength{\unitlength}{0.8pt}
 \begin{picture}(411.70871,454.10507)
 \put(0,0){\includegraphics{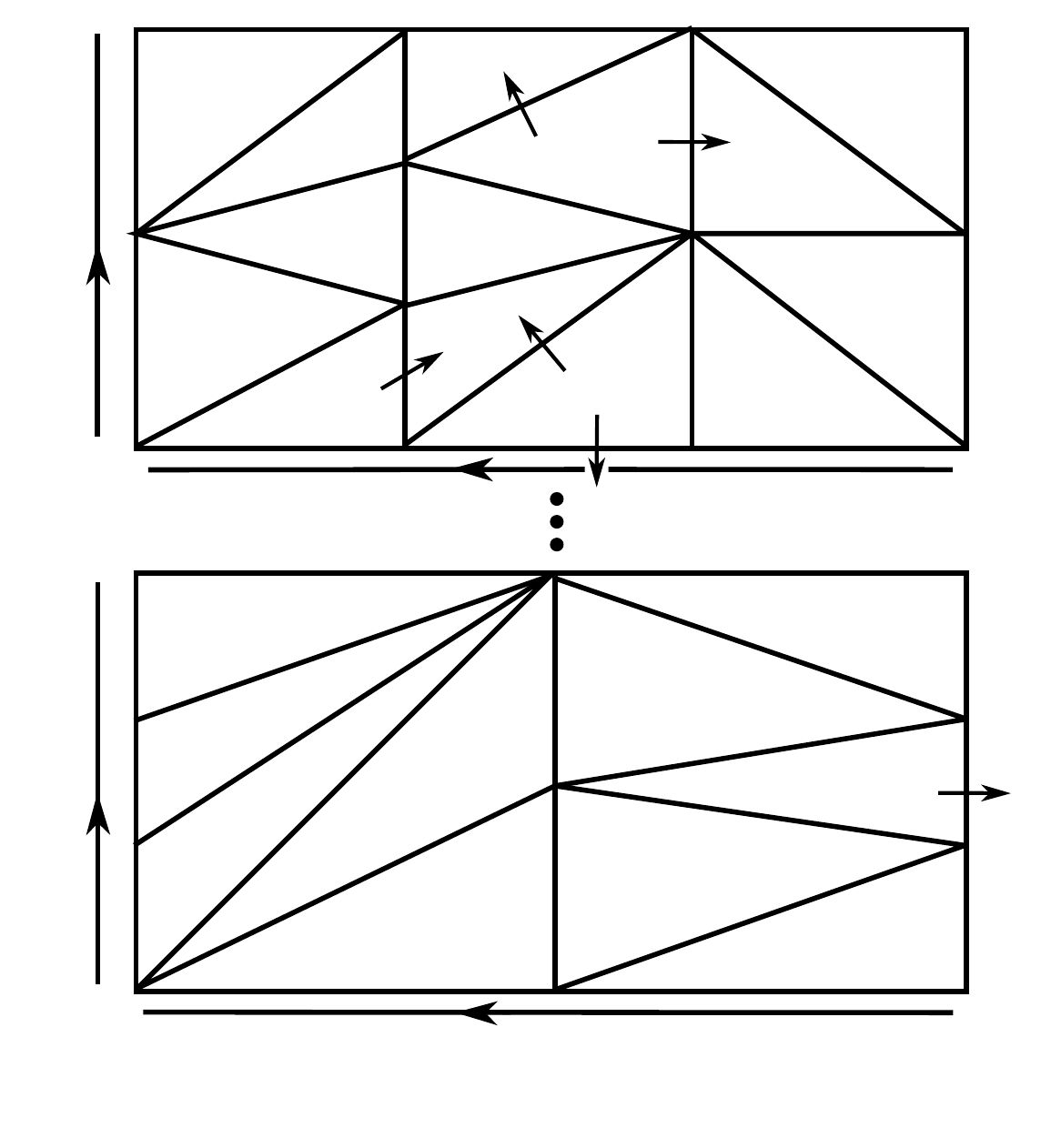}}

\definecolor{inkcol1}{rgb}{0.0,0.0,0.0}
   \put(1.7130675,361.271896){\rotatebox{360.0}{\makebox(0,0)[tl]{\strut{}{
    \begin{minipage}[h]{127.655896pt}
\textcolor{inkcol1}{\LARGE{$\mu_k^\prime$}}\\
\end{minipage}}}}}%

\definecolor{inkcol1}{rgb}{0.0,0.0,0.0}
   \put(176.631891,260.569136){\rotatebox{360.0}{\makebox(0,0)[tl]{\strut{}{
    \begin{minipage}[h]{127.655896pt}
\textcolor{inkcol1}{\LARGE{$\lambda_k^\prime$}}\\
\end{minipage}}}}}%

\definecolor{inkcol1}{rgb}{0.0,0.0,0.0}
   \put(178.316491,44.153536){\rotatebox{360.0}{\makebox(0,0)[tl]{\strut{}{
    \begin{minipage}[h]{127.655896pt}
\textcolor{inkcol1}{\LARGE{$\lambda_l^\prime$}}\\
\end{minipage}}}}}%

\definecolor{inkcol1}{rgb}{0.0,0.0,0.0}
   \put(2.3805856,144.401576){\rotatebox{360.0}{\makebox(0,0)[tl]{\strut{}{
    \begin{minipage}[h]{127.655896pt}
\textcolor{inkcol1}{\LARGE{$\mu_l^\prime$}}\\
\end{minipage}}}}}%

\definecolor{inkcol1}{rgb}{0.0,0.0,0.0}
   \put(297.325641,400.328296){\rotatebox{360.0}{\makebox(0,0)[tl]{\strut{}{
    \begin{minipage}[h]{127.655896pt}
\textcolor{inkcol1}{\LARGE{$0$}}\\
\end{minipage}}}}}%

\definecolor{inkcol1}{rgb}{0.0,0.0,0.0}
   \put(238.484241,414.470426){\rotatebox{360.0}{\makebox(0,0)[tl]{\strut{}{
    \begin{minipage}[h]{127.655896pt}
\textcolor{inkcol1}{\LARGE{$0$}}\\
\end{minipage}}}}}%

\definecolor{inkcol1}{rgb}{0.0,0.0,0.0}
   \put(176.864941,439.219176){\rotatebox{360.0}{\makebox(0,0)[tl]{\strut{}{
    \begin{minipage}[h]{127.655896pt}
\textcolor{inkcol1}{\LARGE{$2$}}\\
\end{minipage}}}}}%

\definecolor{inkcol1}{rgb}{0.0,0.0,0.0}
   \put(124.842081,310.929786){\rotatebox{360.0}{\makebox(0,0)[tl]{\strut{}{
    \begin{minipage}[h]{127.655896pt}
\textcolor{inkcol1}{\LARGE{$2$}}\\
\end{minipage}}}}}%

\definecolor{inkcol1}{rgb}{0.0,0.0,0.0}
   \put(187.471541,339.214056){\rotatebox{360.0}{\makebox(0,0)[tl]{\strut{}{
    \begin{minipage}[h]{127.655896pt}
\textcolor{inkcol1}{\LARGE{$1$}}\\
\end{minipage}}}}}%

\definecolor{inkcol1}{rgb}{0.0,0.0,0.0}
   \put(239.999471,319.768626){\rotatebox{360.0}{\makebox(0,0)[tl]{\strut{}{
    \begin{minipage}[h]{127.655896pt}
\textcolor{inkcol1}{\LARGE{$1$}}\\
\end{minipage}}}}}%

\definecolor{inkcol1}{rgb}{0.0,0.0,0.0}
   \put(297.578171,153.598536){\rotatebox{360.0}{\makebox(0,0)[tl]{\strut{}{
    \begin{minipage}[h]{127.655896pt}
\textcolor{inkcol1}{\LARGE{$3$}}\\
\end{minipage}}}}}%

\definecolor{inkcol1}{rgb}{0.0,0.0,0.0}
   \put(74.839531,173.549036){\rotatebox{360.0}{\makebox(0,0)[tl]{\strut{}{
    \begin{minipage}[h]{127.655896pt}
\textcolor{inkcol1}{\LARGE{$3$}}\\
\end{minipage}}}}}%

\definecolor{inkcol1}{rgb}{0.0,0.0,0.0}
   \put(361.470311,164.962746){\rotatebox{360.0}{\makebox(0,0)[tl]{\strut{}{
    \begin{minipage}[h]{127.655896pt}
\textcolor{inkcol1}{\LARGE{$b$}}\\
\end{minipage}}}}}%

\definecolor{inkcol1}{rgb}{0.0,0.0,0.0}
   \put(141.257051,324.061766){\rotatebox{360.0}{\makebox(0,0)[tl]{\strut{}{
    \begin{minipage}[h]{127.655896pt}
\textcolor{inkcol1}{\LARGE{$a$}}\\
\end{minipage}}}}}%

\definecolor{inkcol1}{rgb}{0.0,0.0,0.0}
   \put(222.069261,326.334616){\rotatebox{360.0}{\makebox(0,0)[tl]{\strut{}{
    \begin{minipage}[h]{127.655896pt}
\textcolor{inkcol1}{\LARGE{$b$}}\\
\end{minipage}}}}}%

\definecolor{inkcol1}{rgb}{0.0,0.0,0.0}
   \put(238.989321,298.302876){\rotatebox{360.0}{\makebox(0,0)[tl]{\strut{}{
    \begin{minipage}[h]{127.655896pt}
\textcolor{inkcol1}{\LARGE{$a$}}\\
\end{minipage}}}}}%

\definecolor{inkcol1}{rgb}{0.0,0.0,0.0}
   \put(212.977891,417.753416){\rotatebox{360.0}{\makebox(0,0)[tl]{\strut{}{
    \begin{minipage}[h]{127.655896pt}
\textcolor{inkcol1}{\LARGE{$b$}}\\
\end{minipage}}}}}%

\definecolor{inkcol1}{rgb}{0.0,0.0,0.0}
   \put(253.636521,395.024986){\rotatebox{360.0}{\makebox(0,0)[tl]{\strut{}{
    \begin{minipage}[h]{127.655896pt}
\textcolor{inkcol1}{\LARGE{$a$}}\\
\end{minipage}}}}}%

 \end{picture}
\endgroup

%% file: figures_gen/ExtIdenRel_n=2.tex
\begingroup
 \setlength{\unitlength}{0.8pt}
 \begin{picture}(283.84589,136.37314)
 \put(0,0){\includegraphics{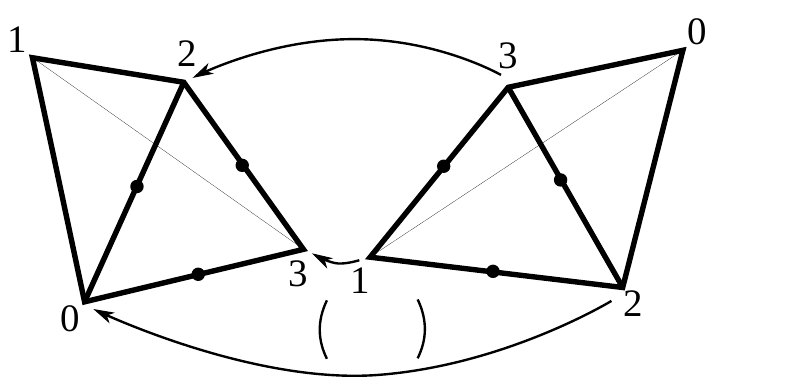}}

\definecolor{inkcol1}{rgb}{0.0,0.0,0.0}
   \put(-0.625263,48.814156){\rotatebox{360.0}{\makebox(0,0)[tl]{\strut{}{
    \begin{minipage}[h]{127.655896pt}
\textcolor{inkcol1}{\large{$\Delta_0$}}\\
\end{minipage}}}}}%

\definecolor{inkcol1}{rgb}{0.0,0.0,0.0}
   \put(233.400307,62.374476){\rotatebox{360.0}{\makebox(0,0)[tl]{\strut{}{
    \begin{minipage}[h]{127.655896pt}
\textcolor{inkcol1}{\large{$\Delta_1$}}\\
\end{minipage}}}}}%

\definecolor{inkcol1}{rgb}{0.0,0.0,0.0}
   \put(115.808107,28.362386){\rotatebox{360.0}{\makebox(0,0)[tl]{\strut{}{
    \begin{minipage}[h]{127.655896pt}
\textcolor{inkcol1}{\normalsize{$0123$}}\\
\end{minipage}}}}}%

\definecolor{inkcol1}{rgb}{0.0,0.0,0.0}
   \put(115.808107,17.755786){\rotatebox{360.0}{\makebox(0,0)[tl]{\strut{}{
    \begin{minipage}[h]{127.655896pt}
\textcolor{inkcol1}{\normalsize{$1302$}}\\
\end{minipage}}}}}%

\definecolor{inkcol1}{rgb}{0.0,0.0,0.0}
   \put(84.294187,24.265576){\rotatebox{360.0}{\makebox(0,0)[tl]{\strut{}{
    \begin{minipage}[h]{127.655896pt}
\textcolor{inkcol1}{\normalsize{$\sigma\!=$}}\\
\end{minipage}}}}}%

\definecolor{inkcol1}{rgb}{0.0,0.0,0.0}
   \put(79.893687,139.382486){\rotatebox{360.0}{\makebox(0,0)[tl]{\strut{}{
    \begin{minipage}[h]{127.655896pt}
\textcolor{inkcol1}{\normalsize{$D_{\lambda^\prime}^{-1}$}}\\
\end{minipage}}}}}%

\definecolor{inkcol1}{rgb}{0.0,0.0,0.0}
   \put(108.921477,62.228276){\rotatebox{360.0}{\makebox(0,0)[tl]{\strut{}{
    \begin{minipage}[h]{127.655896pt}
\textcolor{inkcol1}{\normalsize{$D_{\mu^\prime}$}}\\
\end{minipage}}}}}%

\definecolor{inkcol1}{rgb}{0.0,0.0,0.0}
   \put(42.503947,13.235906){\rotatebox{360.0}{\makebox(0,0)[tl]{\strut{}{
    \begin{minipage}[h]{127.655896pt}
\textcolor{inkcol1}{\normalsize{$D_{\lambda^\prime}$}}\\
\end{minipage}}}}}%

\definecolor{inkcol1}{rgb}{0.0,0.0,0.0}
   \put(87.669867,90.207556){\rotatebox{360.0}{\makebox(0,0)[tl]{\strut{}{
    \begin{minipage}[h]{127.655896pt}
\textcolor{inkcol1}{\normalsize{$(l^\prime)^{-1}m^{\prime}$}}\\
\end{minipage}}}}}%

\definecolor{inkcol1}{rgb}{0.0,0.0,0.0}
   \put(53.533667,55.058646){\rotatebox{360.0}{\makebox(0,0)[tl]{\strut{}{
    \begin{minipage}[h]{127.655896pt}
\textcolor{inkcol1}{\normalsize{$m^{\prime}l^\prime$}}\\
\end{minipage}}}}}%

\definecolor{inkcol1}{rgb}{0.0,0.0,0.0}
   \put(36.613597,83.595456){\rotatebox{360.0}{\makebox(0,0)[tl]{\strut{}{
    \begin{minipage}[h]{127.655896pt}
\textcolor{inkcol1}{\normalsize{$1$}}\\
\end{minipage}}}}}%

 \end{picture}
\endgroup

%% file: figures_gen/ExtIdenRel_n=3.tex
\begingroup
 \setlength{\unitlength}{0.8pt}
 \begin{picture}(274.00656,134.10031)
 \put(0,0){\includegraphics{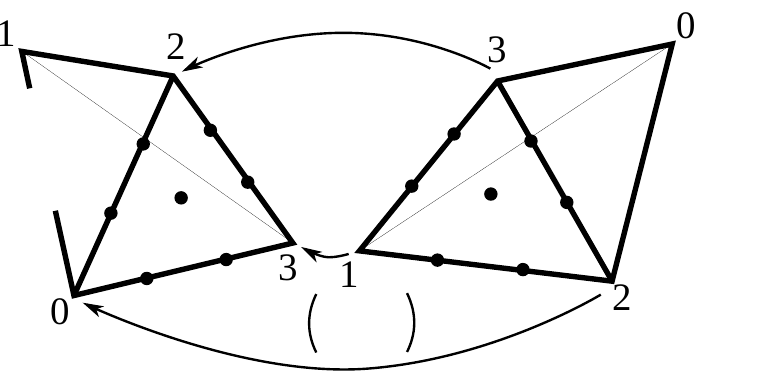}}

\definecolor{inkcol1}{rgb}{0.0,0.0,0.0}
   \put(-4.453141,48.814156){\rotatebox{360.0}{\makebox(0,0)[tl]{\strut{}{
    \begin{minipage}[h]{127.655896pt}
\textcolor{inkcol1}{\large{$\Delta_0$}}\\
\end{minipage}}}}}%

\definecolor{inkcol1}{rgb}{0.0,0.0,0.0}
   \put(229.572429,62.374486){\rotatebox{360.0}{\makebox(0,0)[tl]{\strut{}{
    \begin{minipage}[h]{127.655896pt}
\textcolor{inkcol1}{\large{$\Delta_1$}}\\
\end{minipage}}}}}%

\definecolor{inkcol1}{rgb}{0.0,0.0,0.0}
   \put(111.980229,28.362386){\rotatebox{360.0}{\makebox(0,0)[tl]{\strut{}{
    \begin{minipage}[h]{127.655896pt}
\textcolor{inkcol1}{\normalsize{$0123$}}\\
\end{minipage}}}}}%

\definecolor{inkcol1}{rgb}{0.0,0.0,0.0}
   \put(111.980229,17.755786){\rotatebox{360.0}{\makebox(0,0)[tl]{\strut{}{
    \begin{minipage}[h]{127.655896pt}
\textcolor{inkcol1}{\normalsize{$1302$}}\\
\end{minipage}}}}}%

\definecolor{inkcol1}{rgb}{0.0,0.0,0.0}
   \put(80.466309,24.265576){\rotatebox{360.0}{\makebox(0,0)[tl]{\strut{}{
    \begin{minipage}[h]{127.655896pt}
\textcolor{inkcol1}{\normalsize{$\sigma\!=$}}\\
\end{minipage}}}}}%

\definecolor{inkcol1}{rgb}{0.0,0.0,0.0}
   \put(78.843729,137.109656){\rotatebox{360.0}{\makebox(0,0)[tl]{\strut{}{
    \begin{minipage}[h]{127.655896pt}
\textcolor{inkcol1}{\normalsize{$D_{\mu^\prime}$}}\\
\end{minipage}}}}}%

\definecolor{inkcol1}{rgb}{0.0,0.0,0.0}
   \put(112.922279,59.702906){\rotatebox{360.0}{\makebox(0,0)[tl]{\strut{}{
    \begin{minipage}[h]{127.655896pt}
\textcolor{inkcol1}{\normalsize{$I$}}\\
\end{minipage}}}}}%

\definecolor{inkcol1}{rgb}{0.0,0.0,0.0}
   \put(38.67607,13.235906){\rotatebox{360.0}{\makebox(0,0)[tl]{\strut{}{
    \begin{minipage}[h]{127.655896pt}
\textcolor{inkcol1}{\normalsize{$D_{\lambda^\prime}$}}\\
\end{minipage}}}}}%

\definecolor{inkcol1}{rgb}{0.0,0.0,0.0}
   \put(75.734579,101.309196){\rotatebox{360.0}{\makebox(0,0)[tl]{\strut{}{
    \begin{minipage}[h]{127.655896pt}
\textcolor{inkcol1}{\footnotesize{$m_1^\prime m_2^\prime$}}\\
\end{minipage}}}}}%

\definecolor{inkcol1}{rgb}{0.0,0.0,0.0}
   \put(5.89397,99.561376){\rotatebox{360.0}{\makebox(0,0)[tl]{\strut{}{
    \begin{minipage}[h]{127.655896pt}
\textcolor{inkcol1}{\footnotesize{$m_1^\prime m_2^\prime l_1^\prime $}}\\
\end{minipage}}}}}%

\definecolor{inkcol1}{rgb}{0.0,0.0,0.0}
   \put(1.9444,74.392946){\rotatebox{360.0}{\makebox(0,0)[tl]{\strut{}{
    \begin{minipage}[h]{127.655896pt}
\textcolor{inkcol1}{\footnotesize{$m_1^\prime l_1^\prime l_2^\prime$}}\\
\end{minipage}}}}}%

\definecolor{inkcol1}{rgb}{0.0,0.0,0.0}
   \put(76.943769,38.739556){\rotatebox{360.0}{\makebox(0,0)[tl]{\strut{}{
    \begin{minipage}[h]{127.655896pt}
\textcolor{inkcol1}{\footnotesize{$l_1^\prime$}}\\
\end{minipage}}}}}%

\definecolor{inkcol1}{rgb}{0.0,0.0,0.0}
   \put(48.912032,32.678636){\rotatebox{360.0}{\makebox(0,0)[tl]{\strut{}{
    \begin{minipage}[h]{127.655896pt}
\textcolor{inkcol1}{\footnotesize{$l_1^\prime l_2^\prime$}}\\
\end{minipage}}}}}%

\definecolor{inkcol1}{rgb}{0.0,0.0,0.0}
   \put(56.740713,62.730676){\rotatebox{360.0}{\makebox(0,0)[tl]{\strut{}{
    \begin{minipage}[h]{127.655896pt}
\textcolor{inkcol1}{\footnotesize{$m_1^\prime l_1^\prime$}}\\
\end{minipage}}}}}%

\definecolor{inkcol1}{rgb}{0.0,0.0,0.0}
   \put(91.085899,79.650736){\rotatebox{360.0}{\makebox(0,0)[tl]{\strut{}{
    \begin{minipage}[h]{127.655896pt}
\textcolor{inkcol1}{\footnotesize{$m_1^\prime$}}\\
\end{minipage}}}}}%

 \end{picture}
\endgroup

%% file: figures_gen/m004Triangulation.tex
\begingroup
 \setlength{\unitlength}{0.8pt}
 \begin{picture}(627.5675,306.49551)
 \put(0,0){\includegraphics{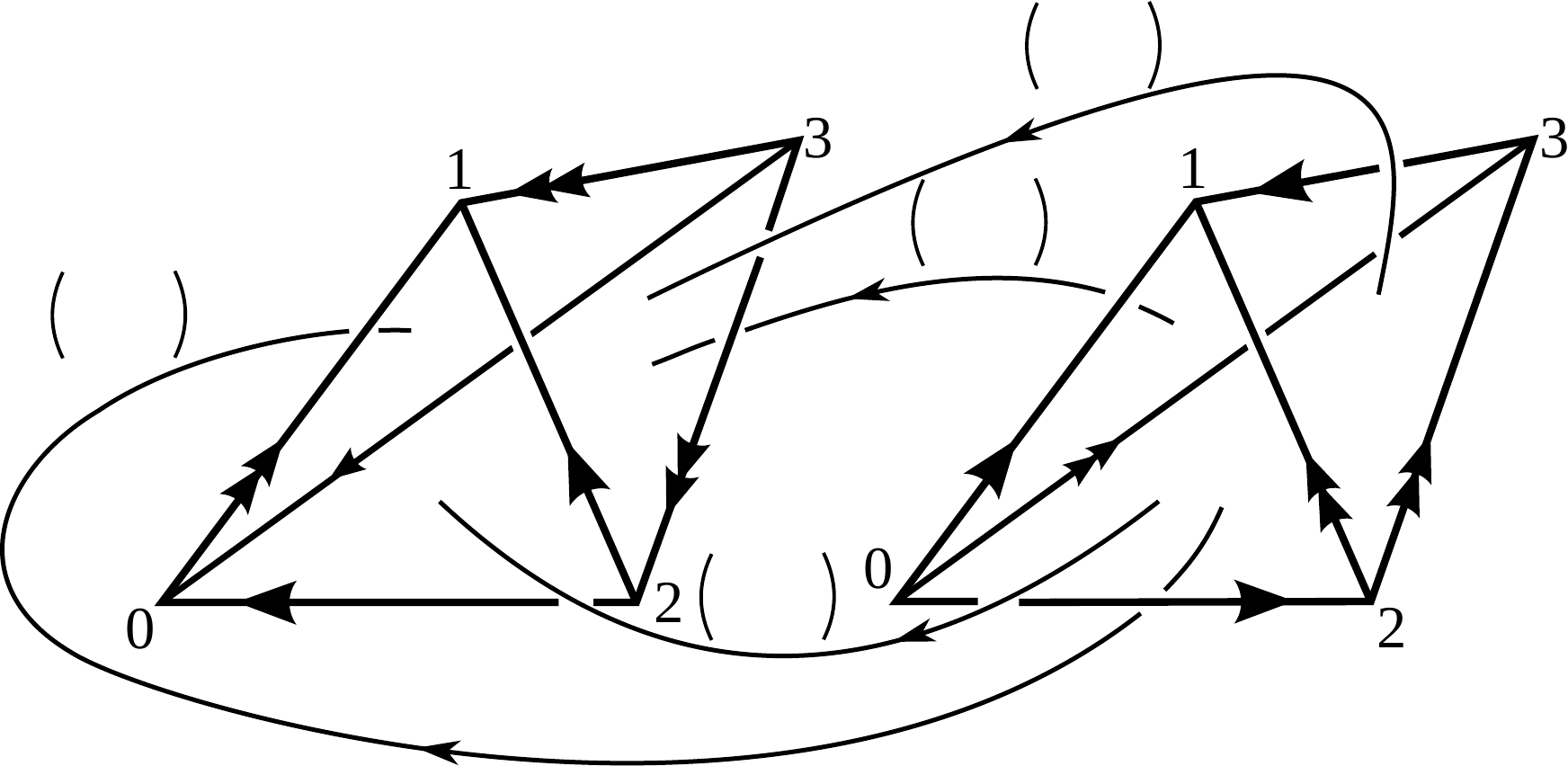}}

\definecolor{inkcol1}{rgb}{0.0,0.0,0.0}
   \put(369.309157,231.294616){\rotatebox{360.0}{\makebox(0,0)[tl]{\strut{}{
    \begin{minipage}[h]{127.655896pt}
\textcolor{inkcol1}{\Large{$0123$}}\\
\end{minipage}}}}}%

\definecolor{inkcol1}{rgb}{0.0,0.0,0.0}
   \put(369.309157,215.637246){\rotatebox{360.0}{\makebox(0,0)[tl]{\strut{}{
    \begin{minipage}[h]{127.655896pt}
\textcolor{inkcol1}{\Large{$3012$}}\\
\end{minipage}}}}}%

\definecolor{inkcol1}{rgb}{0.0,0.0,0.0}
   \put(414.911957,302.093306){\rotatebox{360.0}{\makebox(0,0)[tl]{\strut{}{
    \begin{minipage}[h]{127.655896pt}
\textcolor{inkcol1}{\Large{$0123$}}\\
\end{minipage}}}}}%

\definecolor{inkcol1}{rgb}{0.0,0.0,0.0}
   \put(414.911957,286.435936){\rotatebox{360.0}{\makebox(0,0)[tl]{\strut{}{
    \begin{minipage}[h]{127.655896pt}
\textcolor{inkcol1}{\Large{$0132$}}\\
\end{minipage}}}}}%

\definecolor{inkcol1}{rgb}{0.0,0.0,0.0}
   \put(284.515757,81.407926){\rotatebox{360.0}{\makebox(0,0)[tl]{\strut{}{
    \begin{minipage}[h]{127.655896pt}
\textcolor{inkcol1}{\Large{$0123$}}\\
\end{minipage}}}}}%

\definecolor{inkcol1}{rgb}{0.0,0.0,0.0}
   \put(284.515757,65.750556){\rotatebox{360.0}{\makebox(0,0)[tl]{\strut{}{
    \begin{minipage}[h]{127.655896pt}
\textcolor{inkcol1}{\Large{$2103$}}\\
\end{minipage}}}}}%

\definecolor{inkcol1}{rgb}{0.0,0.0,0.0}
   \put(383.122777,190.427076){\rotatebox{360.0}{\makebox(0,0)[tl]{\strut{}{
    \begin{minipage}[h]{127.655896pt}
\textcolor{inkcol1}{\LARGE{$a$}}\\
\end{minipage}}}}}%

\definecolor{inkcol1}{rgb}{0.0,0.0,0.0}
   \put(437.408497,261.855646){\rotatebox{360.0}{\makebox(0,0)[tl]{\strut{}{
    \begin{minipage}[h]{127.655896pt}
\textcolor{inkcol1}{\LARGE{$d$}}\\
\end{minipage}}}}}%

\definecolor{inkcol1}{rgb}{0.0,0.0,0.0}
   \put(298.837077,39.355646){\rotatebox{360.0}{\makebox(0,0)[tl]{\strut{}{
    \begin{minipage}[h]{127.655896pt}
\textcolor{inkcol1}{\LARGE{$c$}}\\
\end{minipage}}}}}%

\definecolor{inkcol1}{rgb}{0.0,0.0,0.0}
   \put(24.747737,194.283086){\rotatebox{360.0}{\makebox(0,0)[tl]{\strut{}{
    \begin{minipage}[h]{127.655896pt}
\textcolor{inkcol1}{\Large{$0123$}}\\
\end{minipage}}}}}%

\definecolor{inkcol1}{rgb}{0.0,0.0,0.0}
   \put(24.747737,178.625716){\rotatebox{360.0}{\makebox(0,0)[tl]{\strut{}{
    \begin{minipage}[h]{127.655896pt}
\textcolor{inkcol1}{\Large{$3201$}}\\
\end{minipage}}}}}%

\definecolor{inkcol1}{rgb}{0.0,0.0,0.0}
   \put(49.408496,148.998496){\rotatebox{360.0}{\makebox(0,0)[tl]{\strut{}{
    \begin{minipage}[h]{127.655896pt}
\textcolor{inkcol1}{\LARGE{$b$}}\\
\end{minipage}}}}}%

 \end{picture}
\endgroup

%% file: figures_gen/m004Cusp.tex
\begingroup
 \setlength{\unitlength}{0.8pt}
 \begin{picture}(534.47388,230.1161)
 \put(0,0){\includegraphics{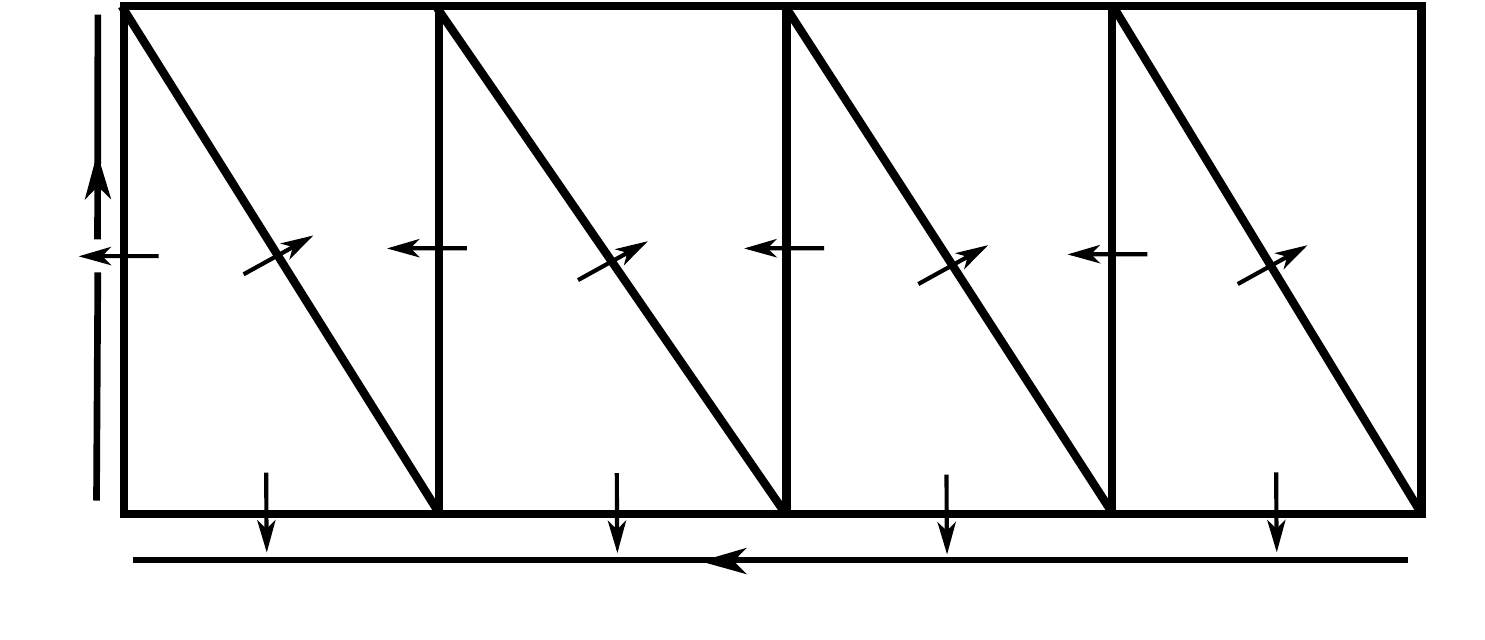}}

\definecolor{inkcol1}{rgb}{0.0,0.0,0.0}
   \put(-0.162374,177.847106){\rotatebox{360.0}{\makebox(0,0)[tl]{\strut{}{
    \begin{minipage}[h]{127.655896pt}
\textcolor{inkcol1}{\huge{$\mu^\prime$}}\\
\end{minipage}}}}}%

\definecolor{inkcol1}{rgb}{0.0,0.0,0.0}
   \put(254.013182,23.450416){\rotatebox{360.0}{\makebox(0,0)[tl]{\strut{}{
    \begin{minipage}[h]{127.655896pt}
\textcolor{inkcol1}{\huge{$\lambda^\prime$}}\\
\end{minipage}}}}}%

\definecolor{inkcol1}{rgb}{0.0,0.0,0.0}
   \put(49.518816,137.710596){\rotatebox{360.0}{\makebox(0,0)[tl]{\strut{}{
    \begin{minipage}[h]{127.655896pt}
\textcolor{inkcol1}{\huge{$a$}}\\
\end{minipage}}}}}%

\definecolor{inkcol1}{rgb}{0.0,0.0,0.0}
   \put(81.473272,131.745236){\rotatebox{360.0}{\makebox(0,0)[tl]{\strut{}{
    \begin{minipage}[h]{127.655896pt}
\textcolor{inkcol1}{\huge{$d$}}\\
\end{minipage}}}}}%

\definecolor{inkcol1}{rgb}{0.0,0.0,0.0}
   \put(162.942452,142.199546){\rotatebox{360.0}{\makebox(0,0)[tl]{\strut{}{
    \begin{minipage}[h]{127.655896pt}
\textcolor{inkcol1}{\huge{$b$}}\\
\end{minipage}}}}}%

\definecolor{inkcol1}{rgb}{0.0,0.0,0.0}
   \put(200.748832,132.127756){\rotatebox{360.0}{\makebox(0,0)[tl]{\strut{}{
    \begin{minipage}[h]{127.655896pt}
\textcolor{inkcol1}{\huge{$d$}}\\
\end{minipage}}}}}%

\definecolor{inkcol1}{rgb}{0.0,0.0,0.0}
   \put(75.131362,69.050886){\rotatebox{360.0}{\makebox(0,0)[tl]{\strut{}{
    \begin{minipage}[h]{127.655896pt}
\textcolor{inkcol1}{\huge{$c$}}\\
\end{minipage}}}}}%

\definecolor{inkcol1}{rgb}{0.0,0.0,0.0}
   \put(201.202782,68.336596){\rotatebox{360.0}{\makebox(0,0)[tl]{\strut{}{
    \begin{minipage}[h]{127.655896pt}
\textcolor{inkcol1}{\huge{$a$}}\\
\end{minipage}}}}}%

\definecolor{inkcol1}{rgb}{0.0,0.0,0.0}
   \put(319.799602,67.979456){\rotatebox{360.0}{\makebox(0,0)[tl]{\strut{}{
    \begin{minipage}[h]{127.655896pt}
\textcolor{inkcol1}{\huge{$a$}}\\
\end{minipage}}}}}%

\definecolor{inkcol1}{rgb}{0.0,0.0,0.0}
   \put(438.345642,67.979456){\rotatebox{360.0}{\makebox(0,0)[tl]{\strut{}{
    \begin{minipage}[h]{127.655896pt}
\textcolor{inkcol1}{\huge{$c$}}\\
\end{minipage}}}}}%

\definecolor{inkcol1}{rgb}{0.0,0.0,0.0}
   \put(288.840212,138.945006){\rotatebox{360.0}{\makebox(0,0)[tl]{\strut{}{
    \begin{minipage}[h]{127.655896pt}
\textcolor{inkcol1}{\huge{$c$}}\\
\end{minipage}}}}}%

\definecolor{inkcol1}{rgb}{0.0,0.0,0.0}
   \put(330.258622,130.280766){\rotatebox{360.0}{\makebox(0,0)[tl]{\strut{}{
    \begin{minipage}[h]{127.655896pt}
\textcolor{inkcol1}{\huge{$b$}}\\
\end{minipage}}}}}%

\definecolor{inkcol1}{rgb}{0.0,0.0,0.0}
   \put(445.488492,128.826816){\rotatebox{360.0}{\makebox(0,0)[tl]{\strut{}{
    \begin{minipage}[h]{127.655896pt}
\textcolor{inkcol1}{\huge{$b$}}\\
\end{minipage}}}}}%

\definecolor{inkcol1}{rgb}{0.0,0.0,0.0}
   \put(405.728922,141.650296){\rotatebox{360.0}{\makebox(0,0)[tl]{\strut{}{
    \begin{minipage}[h]{127.655896pt}
\textcolor{inkcol1}{\huge{$d$}}\\
\end{minipage}}}}}%

\definecolor{inkcol1}{rgb}{0.0,0.0,0.0}
   \put(70.863352,106.296946){\rotatebox{360.0}{\makebox(0,0)[tl]{\strut{}{
    \begin{minipage}[h]{127.655896pt}
\textcolor{inkcol1}{\huge{$1$}}\\
\end{minipage}}}}}%

\definecolor{inkcol1}{rgb}{0.0,0.0,0.0}
   \put(113.289762,185.086866){\rotatebox{360.0}{\makebox(0,0)[tl]{\strut{}{
    \begin{minipage}[h]{127.655896pt}
\textcolor{inkcol1}{\huge{$1$}}\\
\end{minipage}}}}}%

\definecolor{inkcol1}{rgb}{0.0,0.0,0.0}
   \put(189.051202,106.296946){\rotatebox{360.0}{\makebox(0,0)[tl]{\strut{}{
    \begin{minipage}[h]{127.655896pt}
\textcolor{inkcol1}{\huge{$3$}}\\
\end{minipage}}}}}%

\definecolor{inkcol1}{rgb}{0.0,0.0,0.0}
   \put(226.931922,186.602086){\rotatebox{360.0}{\makebox(0,0)[tl]{\strut{}{
    \begin{minipage}[h]{127.655896pt}
\textcolor{inkcol1}{\huge{$2$}}\\
\end{minipage}}}}}%

\definecolor{inkcol1}{rgb}{0.0,0.0,0.0}
   \put(307.239052,106.296946){\rotatebox{360.0}{\makebox(0,0)[tl]{\strut{}{
    \begin{minipage}[h]{127.655896pt}
\textcolor{inkcol1}{\huge{$0$}}\\
\end{minipage}}}}}%

\definecolor{inkcol1}{rgb}{0.0,0.0,0.0}
   \put(350.170532,187.109146){\rotatebox{360.0}{\makebox(0,0)[tl]{\strut{}{
    \begin{minipage}[h]{127.655896pt}
\textcolor{inkcol1}{\huge{$3$}}\\
\end{minipage}}}}}%

\definecolor{inkcol1}{rgb}{0.0,0.0,0.0}
   \put(423.406582,106.296946){\rotatebox{360.0}{\makebox(0,0)[tl]{\strut{}{
    \begin{minipage}[h]{127.655896pt}
\textcolor{inkcol1}{\huge{$2$}}\\
\end{minipage}}}}}%

\definecolor{inkcol1}{rgb}{0.0,0.0,0.0}
   \put(464.822852,187.109146){\rotatebox{360.0}{\makebox(0,0)[tl]{\strut{}{
    \begin{minipage}[h]{127.655896pt}
\textcolor{inkcol1}{\huge{$0$}}\\
\end{minipage}}}}}%

 \end{picture}
\endgroup

%% file: figures_gen/m004Sl2Extended.tex
\begingroup
 \setlength{\unitlength}{0.8pt}
 \begin{picture}(576.31873,212.79555)
 \put(0,0){\includegraphics{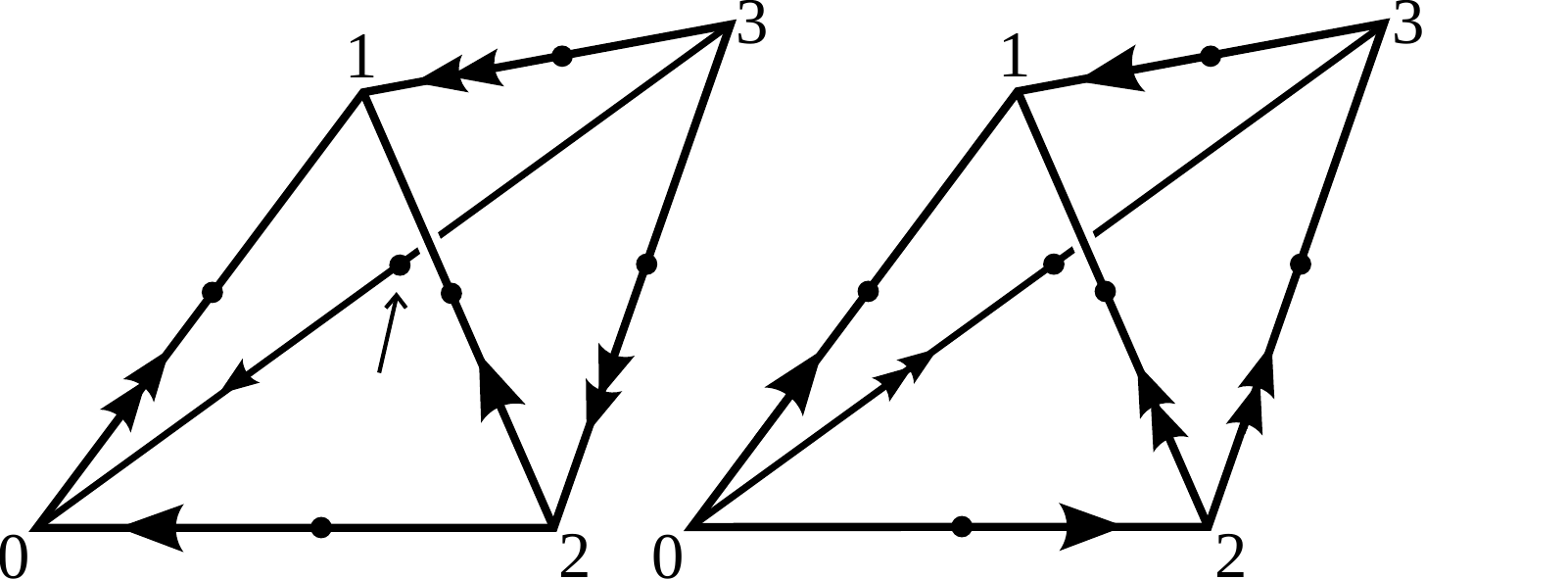}}

\definecolor{inkcol1}{rgb}{0.0,0.0,0.0}
   \put(482.564762,137.555388){\rotatebox{360.0}{\makebox(0,0)[tl]{\strut{}{
    \begin{minipage}[h]{129.560208pt}
\textcolor{inkcol1}{\LARGE{$-m'^2y$}}\\
\end{minipage}}}}}%

\definecolor{inkcol1}{rgb}{0.0,0.0,0.0}
   \put(411.850482,127.198258){\rotatebox{360.0}{\makebox(0,0)[tl]{\strut{}{
    \begin{minipage}[h]{129.560208pt}
\textcolor{inkcol1}{\LARGE{$y$}}\\
\end{minipage}}}}}%

\definecolor{inkcol1}{rgb}{0.0,0.0,0.0}
   \put(352.074032,139.341108){\rotatebox{360.0}{\makebox(0,0)[tl]{\strut{}{
    \begin{minipage}[h]{129.560208pt}
\textcolor{inkcol1}{\LARGE{$-y$}}\\
\end{minipage}}}}}%

\definecolor{inkcol1}{rgb}{0.0,0.0,0.0}
   \put(305.421912,132.555398){\rotatebox{360.0}{\makebox(0,0)[tl]{\strut{}{
    \begin{minipage}[h]{129.560208pt}
\textcolor{inkcol1}{\LARGE{$x$}}\\
\end{minipage}}}}}%

\definecolor{inkcol1}{rgb}{0.0,0.0,0.0}
   \put(307.120392,47.688978){\rotatebox{360.0}{\makebox(0,0)[tl]{\strut{}{
    \begin{minipage}[h]{129.560208pt}
\textcolor{inkcol1}{\LARGE{$m'l'^{-1}x$}}\\
\end{minipage}}}}}%

\definecolor{inkcol1}{rgb}{0.0,0.0,0.0}
   \put(404.350482,221.126828){\rotatebox{360.0}{\makebox(0,0)[tl]{\strut{}{
    \begin{minipage}[h]{129.560208pt}
\textcolor{inkcol1}{\LARGE{$-m'x$}}\\
\end{minipage}}}}}%

\definecolor{inkcol1}{rgb}{0.0,0.0,0.0}
   \put(188.279052,215.769678){\rotatebox{360.0}{\makebox(0,0)[tl]{\strut{}{
    \begin{minipage}[h]{129.560208pt}
\textcolor{inkcol1}{\LARGE{$y$}}\\
\end{minipage}}}}}%

\definecolor{inkcol1}{rgb}{0.0,0.0,0.0}
   \put(160.779052,131.483968){\rotatebox{360.0}{\makebox(0,0)[tl]{\strut{}{
    \begin{minipage}[h]{129.560208pt}
\textcolor{inkcol1}{\LARGE{$-m'x$}}\\
\end{minipage}}}}}%

\definecolor{inkcol1}{rgb}{0.0,0.0,0.0}
   \put(73.489887,74.431088){\rotatebox{360.0}{\makebox(0,0)[tl]{\strut{}{
    \begin{minipage}[h]{129.560208pt}
\textcolor{inkcol1}{\LARGE{$-m'l'^{-1}x$}}\\
\end{minipage}}}}}%

\definecolor{inkcol1}{rgb}{0.0,0.0,0.0}
   \put(5.023877,138.387828){\rotatebox{360.0}{\makebox(0,0)[tl]{\strut{}{
    \begin{minipage}[h]{129.560208pt}
\textcolor{inkcol1}{\LARGE{$-m'^2y$}}\\
\end{minipage}}}}}%

\definecolor{inkcol1}{rgb}{0.0,0.0,0.0}
   \put(68.936997,45.412538){\rotatebox{360.0}{\makebox(0,0)[tl]{\strut{}{
    \begin{minipage}[h]{129.560208pt}
\textcolor{inkcol1}{\LARGE{$-m'^2l'^{-1}x$}}\\
\end{minipage}}}}}%

\definecolor{inkcol1}{rgb}{0.0,0.0,0.0}
   \put(242.636202,134.891608){\rotatebox{360.0}{\makebox(0,0)[tl]{\strut{}{
    \begin{minipage}[h]{129.560208pt}
\textcolor{inkcol1}{\LARGE{$m'^2y$}}\\
\end{minipage}}}}}%

 \end{picture}
\endgroup

%% file: figures_gen/ExtendedPtolemy_n=3.tex
\begingroup
 \setlength{\unitlength}{0.8pt}
 \begin{picture}(573.3963,214.40688)
 \put(0,0){\includegraphics{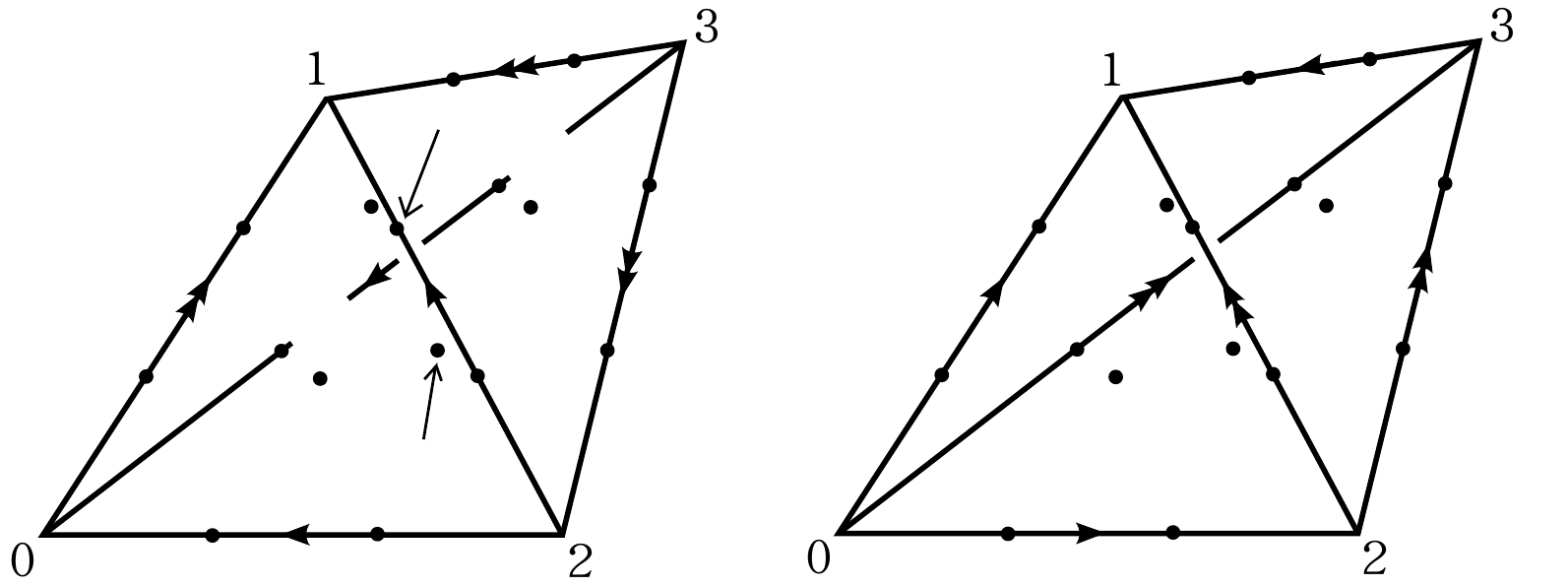}}

\definecolor{inkcol1}{rgb}{0.0,0.0,0.0}
   \put(72.832837,105.562426){\rotatebox{360.0}{\makebox(0,0)[tl]{\strut{}{
    \begin{minipage}[h]{127.655896pt}
\textcolor{inkcol1}{\normalsize{$m_1^3m_2^2l_1l_2x_1$}}\\
\end{minipage}}}}}%

\definecolor{inkcol1}{rgb}{0.0,0.0,0.0}
   \put(167.697697,166.799526){\rotatebox{360.0}{\makebox(0,0)[tl]{\strut{}{
    \begin{minipage}[h]{127.655896pt}
\textcolor{inkcol1}{\normalsize{$m_1^3m_2l_1x_0$}}\\
\end{minipage}}}}}%

\definecolor{inkcol1}{rgb}{0.0,0.0,0.0}
   \put(123.830037,53.922406){\rotatebox{360.0}{\makebox(0,0)[tl]{\strut{}{
    \begin{minipage}[h]{127.655896pt}
\textcolor{inkcol1}{\normalsize{$m_1^4l_1f_2$}}\\
\end{minipage}}}}}%

\definecolor{inkcol1}{rgb}{0.0,0.0,0.0}
   \put(148.042617,200.445676){\rotatebox{360.0}{\makebox(0,0)[tl]{\strut{}{
    \begin{minipage}[h]{127.655896pt}
\textcolor{inkcol1}{\normalsize{$y_1$}}\\
\end{minipage}}}}}%

\definecolor{inkcol1}{rgb}{0.0,0.0,0.0}
   \put(193.320907,206.629146){\rotatebox{360.0}{\makebox(0,0)[tl]{\strut{}{
    \begin{minipage}[h]{127.655896pt}
\textcolor{inkcol1}{\normalsize{$y_0$}}\\
\end{minipage}}}}}%

\definecolor{inkcol1}{rgb}{0.0,0.0,0.0}
   \put(37.757417,153.596146){\rotatebox{360.0}{\makebox(0,0)[tl]{\strut{}{
    \begin{minipage}[h]{127.655896pt}
\textcolor{inkcol1}{\normalsize{$m_1^2m_2y_1$}}\\
\end{minipage}}}}}%

\definecolor{inkcol1}{rgb}{0.0,0.0,0.0}
   \put(-0.123313,98.542816){\rotatebox{360.0}{\makebox(0,0)[tl]{\strut{}{
    \begin{minipage}[h]{127.655896pt}
\textcolor{inkcol1}{\normalsize{$m_1^2m_2y_0$}}\\
\end{minipage}}}}}%

\definecolor{inkcol1}{rgb}{0.0,0.0,0.0}
   \put(35.232037,15.205246){\rotatebox{360.0}{\makebox(0,0)[tl]{\strut{}{
    \begin{minipage}[h]{127.655896pt}
\textcolor{inkcol1}{\normalsize{$m_1^4m_2^3l_1l_2x_1$}}\\
\end{minipage}}}}}%

\definecolor{inkcol1}{rgb}{0.0,0.0,0.0}
   \put(114.529007,14.700166){\rotatebox{360.0}{\makebox(0,0)[tl]{\strut{}{
    \begin{minipage}[h]{127.655896pt}
\textcolor{inkcol1}{\normalsize{$m_1^4m_2l_1x_0$}}\\
\end{minipage}}}}}%

\definecolor{inkcol1}{rgb}{0.0,0.0,0.0}
   \put(115.947057,150.892186){\rotatebox{360.0}{\makebox(0,0)[tl]{\strut{}{
    \begin{minipage}[h]{127.655896pt}
\textcolor{inkcol1}{\normalsize{$f_1$}}\\
\end{minipage}}}}}%

\definecolor{inkcol1}{rgb}{0.0,0.0,0.0}
   \put(135.139977,179.382216){\rotatebox{360.0}{\makebox(0,0)[tl]{\strut{}{
    \begin{minipage}[h]{127.655896pt}
\textcolor{inkcol1}{\normalsize{$m_1m_2x_1$}}\\
\end{minipage}}}}}%

\definecolor{inkcol1}{rgb}{0.0,0.0,0.0}
   \put(241.176247,156.626596){\rotatebox{360.0}{\makebox(0,0)[tl]{\strut{}{
    \begin{minipage}[h]{127.655896pt}
\textcolor{inkcol1}{\normalsize{$m_1^2m_2y_0$}}\\
\end{minipage}}}}}%

\definecolor{inkcol1}{rgb}{0.0,0.0,0.0}
   \put(226.529037,97.532666){\rotatebox{360.0}{\makebox(0,0)[tl]{\strut{}{
    \begin{minipage}[h]{127.655896pt}
\textcolor{inkcol1}{\normalsize{$m_1^2m_2y_1$}}\\
\end{minipage}}}}}%

\definecolor{inkcol1}{rgb}{0.0,0.0,0.0}
   \put(174.174407,89.033036){\rotatebox{360.0}{\makebox(0,0)[tl]{\strut{}{
    \begin{minipage}[h]{127.655896pt}
\textcolor{inkcol1}{\normalsize{$m_1x_0$}}\\
\end{minipage}}}}}%

\definecolor{inkcol1}{rgb}{0.0,0.0,0.0}
   \put(84.097517,70.763636){\rotatebox{360.0}{\makebox(0,0)[tl]{\strut{}{
    \begin{minipage}[h]{127.655896pt}
\textcolor{inkcol1}{\normalsize{$-m_1^2f_3$}}\\
\end{minipage}}}}}%

\definecolor{inkcol1}{rgb}{0.0,0.0,0.0}
   \put(185.178357,139.162266){\rotatebox{360.0}{\makebox(0,0)[tl]{\strut{}{
    \begin{minipage}[h]{127.655896pt}
\textcolor{inkcol1}{\normalsize{$-f_0$}}\\
\end{minipage}}}}}%

\definecolor{inkcol1}{rgb}{0.0,0.0,0.0}
   \put(382.253547,101.862416){\rotatebox{360.0}{\makebox(0,0)[tl]{\strut{}{
    \begin{minipage}[h]{127.655896pt}
\textcolor{inkcol1}{\normalsize{$y_0$}}\\
\end{minipage}}}}}%

\definecolor{inkcol1}{rgb}{0.0,0.0,0.0}
   \put(463.733887,162.341906){\rotatebox{360.0}{\makebox(0,0)[tl]{\strut{}{
    \begin{minipage}[h]{127.655896pt}
\textcolor{inkcol1}{\normalsize{$y_1$}}\\
\end{minipage}}}}}%

\definecolor{inkcol1}{rgb}{0.0,0.0,0.0}
   \put(430.472847,100.224946){\rotatebox{360.0}{\makebox(0,0)[tl]{\strut{}{
    \begin{minipage}[h]{127.655896pt}
\textcolor{inkcol1}{\normalsize{$f_1$}}\\
\end{minipage}}}}}%

\definecolor{inkcol1}{rgb}{0.0,0.0,0.0}
   \put(427.158757,202.554046){\rotatebox{360.0}{\makebox(0,0)[tl]{\strut{}{
    \begin{minipage}[h]{127.655896pt}
\textcolor{inkcol1}{\normalsize{$m_1m_2x_1$}}\\
\end{minipage}}}}}%

\definecolor{inkcol1}{rgb}{0.0,0.0,0.0}
   \put(484.306337,207.727366){\rotatebox{360.0}{\makebox(0,0)[tl]{\strut{}{
    \begin{minipage}[h]{127.655896pt}
\textcolor{inkcol1}{\normalsize{$m_1x_0$}}\\
\end{minipage}}}}}%

\definecolor{inkcol1}{rgb}{0.0,0.0,0.0}
   \put(360.815197,144.087756){\rotatebox{360.0}{\makebox(0,0)[tl]{\strut{}{
    \begin{minipage}[h]{127.655896pt}
\textcolor{inkcol1}{\normalsize{$x_1$}}\\
\end{minipage}}}}}%

\definecolor{inkcol1}{rgb}{0.0,0.0,0.0}
   \put(327.480147,91.559816){\rotatebox{360.0}{\makebox(0,0)[tl]{\strut{}{
    \begin{minipage}[h]{127.655896pt}
\textcolor{inkcol1}{\normalsize{$x_0$}}\\
\end{minipage}}}}}%

\definecolor{inkcol1}{rgb}{0.0,0.0,0.0}
   \put(326.217467,15.798386){\rotatebox{360.0}{\makebox(0,0)[tl]{\strut{}{
    \begin{minipage}[h]{127.655896pt}
\textcolor{inkcol1}{\normalsize{$m_1^3m_2l_1x_0$}}\\
\end{minipage}}}}}%

\definecolor{inkcol1}{rgb}{0.0,0.0,0.0}
   \put(405.514437,15.293306){\rotatebox{360.0}{\makebox(0,0)[tl]{\strut{}{
    \begin{minipage}[h]{127.655896pt}
\textcolor{inkcol1}{\normalsize{$m_1^3m_2^2l_1l_2x_1$}}\\
\end{minipage}}}}}%

\definecolor{inkcol1}{rgb}{0.0,0.0,0.0}
   \put(406.174877,151.485326){\rotatebox{360.0}{\makebox(0,0)[tl]{\strut{}{
    \begin{minipage}[h]{127.655896pt}
\textcolor{inkcol1}{\normalsize{$f_2$}}\\
\end{minipage}}}}}%

\definecolor{inkcol1}{rgb}{0.0,0.0,0.0}
   \put(437.237087,145.377636){\rotatebox{360.0}{\makebox(0,0)[tl]{\strut{}{
    \begin{minipage}[h]{127.655896pt}
\textcolor{inkcol1}{\normalsize{$y_1$}}\\
\end{minipage}}}}}%

\definecolor{inkcol1}{rgb}{0.0,0.0,0.0}
   \put(465.159837,89.626176){\rotatebox{360.0}{\makebox(0,0)[tl]{\strut{}{
    \begin{minipage}[h]{127.655896pt}
\textcolor{inkcol1}{\normalsize{$y_0$}}\\
\end{minipage}}}}}%

\definecolor{inkcol1}{rgb}{0.0,0.0,0.0}
   \put(409.680677,88.024286){\rotatebox{360.0}{\makebox(0,0)[tl]{\strut{}{
    \begin{minipage}[h]{127.655896pt}
\textcolor{inkcol1}{\normalsize{$f_3$}}\\
\end{minipage}}}}}%

\definecolor{inkcol1}{rgb}{0.0,0.0,0.0}
   \put(486.770387,141.523176){\rotatebox{360.0}{\makebox(0,0)[tl]{\strut{}{
    \begin{minipage}[h]{127.655896pt}
\textcolor{inkcol1}{\normalsize{$f_0$}}\\
\end{minipage}}}}}%

\definecolor{inkcol1}{rgb}{0.0,0.0,0.0}
   \put(530.741817,157.297676){\rotatebox{360.0}{\makebox(0,0)[tl]{\strut{}{
    \begin{minipage}[h]{127.655896pt}
\textcolor{inkcol1}{\normalsize{$m_1^2m_2y_1$}}\\
\end{minipage}}}}}%

\definecolor{inkcol1}{rgb}{0.0,0.0,0.0}
   \put(515.262117,96.775046){\rotatebox{360.0}{\makebox(0,0)[tl]{\strut{}{
    \begin{minipage}[h]{127.655896pt}
\textcolor{inkcol1}{\normalsize{$m_1^2m_2y_0$}}\\
\end{minipage}}}}}%

 \end{picture}
\endgroup

%% file: figures_gen/FacePairingGenerator.tex
\begingroup
 \setlength{\unitlength}{0.8pt}
 \begin{picture}(333.4212,227.4606)
 \put(0,0){\includegraphics{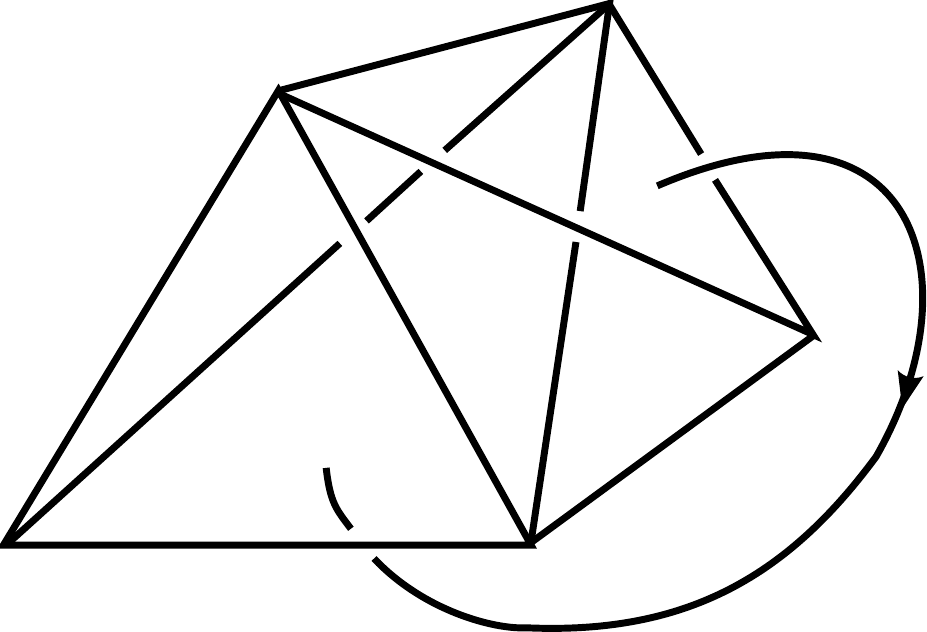}}
 \end{picture}
\endgroup

%% file: figures_gen/FacePairingRelation.tex
\begingroup
 \setlength{\unitlength}{0.8pt}
 \begin{picture}(171.71577,168.68532)
 \put(0,0){\includegraphics{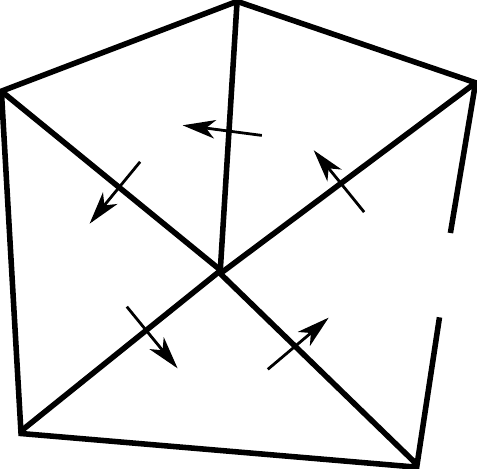}}

\definecolor{inkcol1}{rgb}{0.0,0.0,0.0}
   \put(84.779051,138.463048){\rotatebox{360.0}{\makebox(0,0)[tl]{\strut{}{
    \begin{minipage}[h]{129.560208pt}
\textcolor{inkcol1}{\LARGE{$a$}}\\
\end{minipage}}}}}%

\definecolor{inkcol1}{rgb}{0.0,0.0,0.0}
   \put(33.350481,128.989548){\rotatebox{360.0}{\makebox(0,0)[tl]{\strut{}{
    \begin{minipage}[h]{129.560208pt}
\textcolor{inkcol1}{\LARGE{$b$}}\\
\end{minipage}}}}}%

\definecolor{inkcol1}{rgb}{0.0,0.0,0.0}
   \put(28.707631,61.132408){\rotatebox{360.0}{\makebox(0,0)[tl]{\strut{}{
    \begin{minipage}[h]{129.560208pt}
\textcolor{inkcol1}{\LARGE{$c$}}\\
\end{minipage}}}}}%

\definecolor{inkcol1}{rgb}{0.0,0.0,0.0}
   \put(91.921911,39.346698){\rotatebox{360.0}{\makebox(0,0)[tl]{\strut{}{
    \begin{minipage}[h]{129.560208pt}
\textcolor{inkcol1}{\LARGE{$d$}}\\
\end{minipage}}}}}%

\definecolor{inkcol1}{rgb}{0.0,0.0,0.0}
   \put(111.921911,95.418118){\rotatebox{360.0}{\makebox(0,0)[tl]{\strut{}{
    \begin{minipage}[h]{129.560208pt}
\textcolor{inkcol1}{\LARGE{$e$}}\\
\end{minipage}}}}}%

\definecolor{inkcol1}{rgb}{0.0,0.0,0.0}
   \put(54.064761,34.346688){\rotatebox{360.0}{\makebox(0,0)[tl]{\strut{}{
    \begin{minipage}[h]{129.560208pt}
\textcolor{inkcol1}{\LARGE{$gB$}}\\
\end{minipage}}}}}%

\definecolor{inkcol1}{rgb}{0.0,0.0,0.0}
   \put(92.421911,79.508108){\rotatebox{360.0}{\makebox(0,0)[tl]{\strut{}{
    \begin{minipage}[h]{129.560208pt}
\textcolor{inkcol1}{\LARGE{$\rho(d)gB$}}\\
\end{minipage}}}}}%

 \end{picture}
\endgroup

%% file: figures_gen/FattenedCocycle.tex
\begingroup
 \setlength{\unitlength}{0.8pt}
 \begin{picture}(321.48572,272.20001)
 \put(0,0){\includegraphics{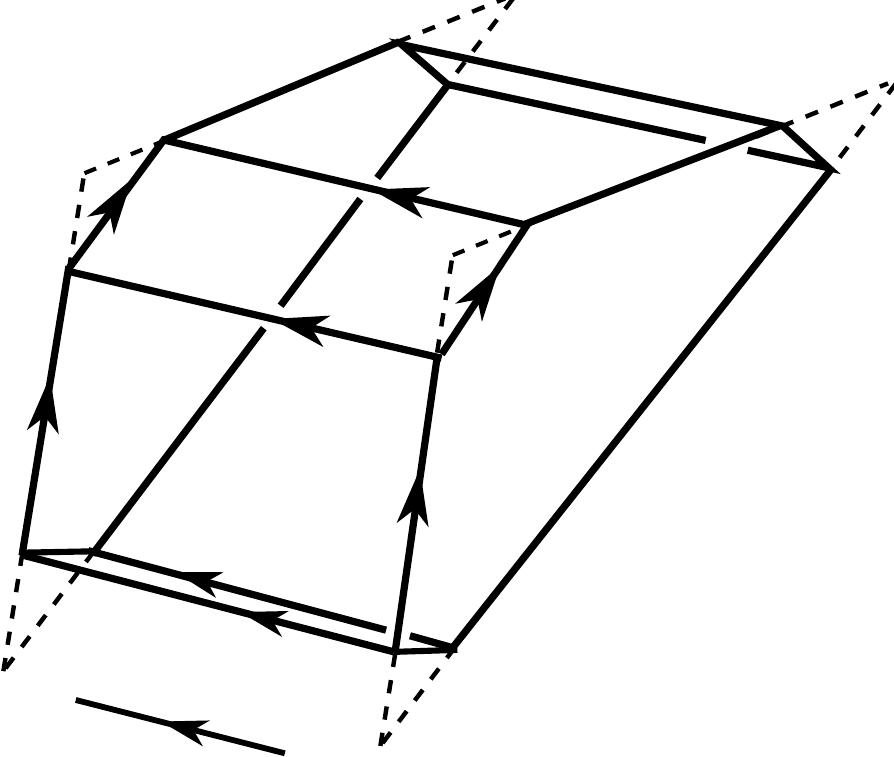}}

\definecolor{inkcol1}{rgb}{0.0,0.0,0.0}
   \put(137.669337,192.988466){\rotatebox{360.0}{\makebox(0,0)[tl]{\strut{}{
    \begin{minipage}[h]{127.655896pt}
\textcolor{inkcol1}{\LARGE{$1$}}\\
\end{minipage}}}}}%

\definecolor{inkcol1}{rgb}{0.0,0.0,0.0}
   \put(2.813977,224.303186){\rotatebox{360.0}{\makebox(0,0)[tl]{\strut{}{
    \begin{minipage}[h]{127.655896pt}
\textcolor{inkcol1}{\LARGE{$1$}}\\
\end{minipage}}}}}%

\definecolor{inkcol1}{rgb}{0.0,0.0,0.0}
   \put(119.286247,4.348496){\rotatebox{360.0}{\makebox(0,0)[tl]{\strut{}{
    \begin{minipage}[h]{127.655896pt}
\textcolor{inkcol1}{\LARGE{$0$}}\\
\end{minipage}}}}}%

\definecolor{inkcol1}{rgb}{0.0,0.0,0.0}
   \put(-15.064043,28.087086){\rotatebox{360.0}{\makebox(0,0)[tl]{\strut{}{
    \begin{minipage}[h]{127.655896pt}
\textcolor{inkcol1}{\LARGE{$0$}}\\
\end{minipage}}}}}%

\definecolor{inkcol1}{rgb}{0.0,0.0,0.0}
   \put(182.925857,296.282576){\rotatebox{360.0}{\makebox(0,0)[tl]{\strut{}{
    \begin{minipage}[h]{127.655896pt}
\textcolor{inkcol1}{\LARGE{$2$}}\\
\end{minipage}}}}}%

\definecolor{inkcol1}{rgb}{0.0,0.0,0.0}
   \put(318.791367,266.988146){\rotatebox{360.0}{\makebox(0,0)[tl]{\strut{}{
    \begin{minipage}[h]{127.655896pt}
\textcolor{inkcol1}{\LARGE{$2$}}\\
\end{minipage}}}}}%

\definecolor{inkcol1}{rgb}{0.0,0.0,0.0}
   \put(176.864937,173.549046){\rotatebox{360.0}{\makebox(0,0)[tl]{\strut{}{
    \begin{minipage}[h]{127.655896pt}
\textcolor{inkcol1}{\LARGE{$\beta_{102,0}$}}\\
\end{minipage}}}}}%

\definecolor{inkcol1}{rgb}{0.0,0.0,0.0}
   \put(156.598537,119.644206){\rotatebox{360.0}{\makebox(0,0)[tl]{\strut{}{
    \begin{minipage}[h]{127.655896pt}
\textcolor{inkcol1}{\LARGE{$\alpha_{01,0}$}}\\
\end{minipage}}}}}%

\definecolor{inkcol1}{rgb}{0.0,0.0,0.0}
   \put(22.799017,149.948776){\rotatebox{360.0}{\makebox(0,0)[tl]{\strut{}{
    \begin{minipage}[h]{127.655896pt}
\textcolor{inkcol1}{\LARGE{$\alpha_{01,1}$}}\\
\end{minipage}}}}}%

\definecolor{inkcol1}{rgb}{0.0,0.0,0.0}
   \put(45.986837,207.022396){\rotatebox{360.0}{\makebox(0,0)[tl]{\strut{}{
    \begin{minipage}[h]{127.655896pt}
\textcolor{inkcol1}{\LARGE{$\beta_{102,1}$}}\\
\end{minipage}}}}}%

\definecolor{inkcol1}{rgb}{0.0,0.0,0.0}
   \put(54.618817,91.359936){\rotatebox{360.0}{\makebox(0,0)[tl]{\strut{}{
    \begin{minipage}[h]{127.655896pt}
\textcolor{inkcol1}{\LARGE{$M_{a0}$}}\\
\end{minipage}}}}}%

\definecolor{inkcol1}{rgb}{0.0,0.0,0.0}
   \put(49.585717,5.358646){\rotatebox{360.0}{\makebox(0,0)[tl]{\strut{}{
    \begin{minipage}[h]{127.655896pt}
\textcolor{inkcol1}{\LARGE{$a$}}\\
\end{minipage}}}}}%

\definecolor{inkcol1}{rgb}{0.0,0.0,0.0}
   \put(87.694867,41.570616){\rotatebox{360.0}{\makebox(0,0)[tl]{\strut{}{
    \begin{minipage}[h]{127.655896pt}
\textcolor{inkcol1}{\LARGE{$M_{a0}$}}\\
\end{minipage}}}}}%

\definecolor{inkcol1}{rgb}{0.0,0.0,0.0}
   \put(104.362387,146.626486){\rotatebox{360.0}{\makebox(0,0)[tl]{\strut{}{
    \begin{minipage}[h]{127.655896pt}
\textcolor{inkcol1}{\LARGE{$M_{a1}$}}\\
\end{minipage}}}}}%

\definecolor{inkcol1}{rgb}{0.0,0.0,0.0}
   \put(151.839557,224.408236){\rotatebox{360.0}{\makebox(0,0)[tl]{\strut{}{
    \begin{minipage}[h]{127.655896pt}
\textcolor{inkcol1}{\LARGE{$M_{a1}$}}\\
\end{minipage}}}}}%

 \end{picture}
\endgroup

%% file: figures_gen/FattenedCuspTriangulation.tex
\begingroup
 \setlength{\unitlength}{0.8pt}
 \begin{picture}(187.71428,178.42857)
 \put(0,0){\includegraphics{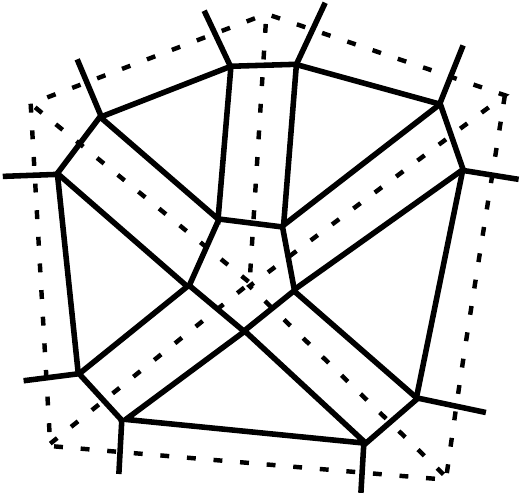}}
 \end{picture}
\endgroup